\documentclass[12pt]{article}%
\usepackage{amsmath}
\usepackage{amsfonts}
\usepackage{amssymb}
\usepackage{graphicx}
\usepackage{array}%
\newtheorem{theorem}{Theorem}

\newtheorem{lemma}[theorem]{Lemma}
\newtheorem{proposition}[theorem]{Proposition}

\newtheorem{remark}[theorem]{Remark}

\newcommand{\epr}{$\Box$}
\newcommand{\abst}{\hspace{1cm}}
\newcommand{\ep}{\varepsilon}
\newcommand{\Om}{\Omega}
\newcommand{\R}{\mathbb{R}}

\makeatletter
\def\@currentlabel{2.1}\label{e:dispaa}
\def\@currentlabel{2.21}\label{e:dispau}
\def\@currentlabel{2.22}\label{e:dispav}
\def\@currentlabel{2.23}\label{e:dispaw}
\def\@currentlabel{2.24}\label{e:dispax}
\def\theequation{\thesection.\@arabic\c@equation}
\makeatother

\def\epsilon{\varepsilon}

\begin{document}

\author{ Henri Berestycki \thanks{\'{E}cole des hautes \'{e}tudes en sciences sociales, CAMS, 190-198, avenue de France, 75244 Paris cedex 13, France. ({\tt hb@ehess.fr}).}
\and
 Juncheng Wei
\thanks{Department of Mathematics,
The Chinese University of Hong Kong,
Shatin, Hong Kong, China, and
Department of Mathematics, University of British Columbia, Vancouver, B.C., Canada, V6T 1Z2. ({\tt wei@math.cuhk.edu.hk}).}
\and
Matthias Winter
\thanks{Brunel University, Department of Mathematical Sciences, Uxbridge UB8 3PH, United Kingdom
({\tt matthias.winter@brunel.ac.uk})}}

\title{Existence of Symmetric and Asymmetric Spikes for a Crime Hotspot Model}
\date{}
\maketitle

\begin{abstract}
We study a crime hotspot model suggested by Short-Bertozzi-Brantingham \cite{sbb}. The aim of this work is to establish rigorously the formation of hotspots in this model representing concentrations of  criminal activity. More precisely,  for the one-dimensional system,
we rigorously prove the existence
of steady states with multiple spikes of the following types:

(i) Multiple spikes of arbitrary number having the same amplitude (symmetric spikes),

(ii) Multiple spikes having different amplitude for the case of one large and one small spike (asymmetric spikes).

We use an approach based on Liapunov-Schmidt reduction and extend it to the
 quasilinear crime hotspot model. Some novel results that allow us to carry out the Liapunov-Schmidt reduction are:
(i) approximation of the quasilinear crime hotspot system on the large scale by the semilinear Schnakenberg model, (ii) estimate of the spatial dependence of the second component on the small scale which is dominated by the quasilinear part of the system.

The paper concludes with an extension to the anisotropic case.

\end{abstract}

\medskip


\medskip


\bigskip

\section{Introduction: The statement of the problem}
\setcounter{equation}{0}

Pattern forming reaction-diffusion systems have been and are applied to many phenomena in the natural sciences. Recent works have also started to use such systems to describe macroscopic social phenomena.
 In this direction, Short, Bertozzi and Brantingham \cite{sbb} have proposed a system of non-linear parabolic partial differential equations to describe the formation of hotspots of criminal activity.
Their equations are  derived from an agent-based lattice model that incorporates the movement of criminals and a given scalar field representing the ``attractiveness of crime''. The system in one dimension reads as follows:
\begin{align}
\nonumber
A_{t}  & =\varepsilon^{2}A_{xx}-A+\rho A+A_{0} (x), \ \mbox{in} \ (-L, L),\\
\rho_{t}  & =D (\rho_{x}-2\frac{\rho}{A}A_{x})_{x}-\rho
A+\gamma (x), \ \mbox{in} \ (-L, L).
\label{sysoriginal}
\end{align}
Here $A$ is the ``attractiveness of crime'' and $\rho$ denotes the density of criminals. The rate at which crimes occur is given by $\rho A$. When this rate increases, the number of criminals is reduced while the attractiveness increases. The second feature is related to the well documented occurrence of repeat offenses. The positive function $A_0 (x)$ is the intrinsic (static) attractiveness which is stationary in time but possibly variable in space. The positive function $\gamma (x)$ is the source term representing the introduction rate of offenders (per unit area). For the precise meanings of the functions $A_0 (x)$ and $\gamma (x)$, we refer to \cite{sbb, sbbt, soptbbc} and the references therein.

This paper is concerned with the mathematical analysis of the one-dimensional version of this system.
Let us describe our approach. Setting
\[
v=\frac{\rho}{A^{2}},
\]
the system is transformed into
\begin{align}
\nonumber
A_{t}  & =\varepsilon^{2}A_{xx}-A+vA^{3}+A_{0} (x) \ \mbox{in} \ (-L, L),\\
(A^2v)_{t}  & =D\left(  A^{2}v_{x}\right)  _{x}-vA^{3}+\gamma (x) \ \mbox{in} \ (-L, L).
\label{sysdyn}
\end{align}
We always consider Neumann boundary conditions
\[
A_x(-L)=A_x(L)=\rho_x(-L)=\rho_x(L)=v_{x}(-L)=v_x(L)=0.
\]
Note that $v$ is well-defined and positive if $A$ and $\rho$ are both positive.

The parameter $0<\ep^2$ represents nearest neighbor interactions in a lattice model for the attractiveness. We assume that it is very small which corresponds to the temporal dependence of attractiveness dominating its spatial dependence. This models the case of attractiveness propagating rather slowly, i.e. much slower than individual criminals. It is a realistic assumption if the criminal spatial profile remains largely unchanged, or, in other words, if the relative crime-intensity does only change very slowly. This appears to be a reasonable assumption since it typically takes decades for dangerous neighborhoods, i.e. those attracting criminals, to evolve into safe ones and vice versa.

Roughly speaking, a $k$ spike solution $(A, v)$ to (\ref{sysdyn}) is such that the component $A$ has exactly $k$ local maximum points.  In this paper, we address the issue of existence of steady states with multiple spikes in the following two cases: Symmetric spikes (same amplitudes) or asymmetric spikes (different amplitudes). Our approach is by rigorous nonlinear analysis. We apply Liapunov-Schmidt reduction to this quasilinear system.

In this approach, to establish the existence of spikes, we derive the following new results:
\begin{itemize}
\item[(i)] Approximation of the crime hotspot system on the large scale of order one by the semi-linear Schnakenberg model (see Section 3, in particular equation (\ref{approx})),
\item[(ii)] Estimate of the spatial dependence of the second component on the small scale of order $\ep$, dominated by the quasilinear part of the system (see Section 6, in particular inequalities  (\ref{estw3}) -- (\ref{estw1})).
\end{itemize}

We remark that asymmetric multiple spike steady states (of $k_1$ small and $k_2$ large spikes) are an intermediate state between two different symmetric multiple spike steady states of $k_1+k_2$ spikes (for which all spikes are fully developed) and $k_2$ spikes (for which the small spikes are gone). These rigorous results shed light on the formation of hotspots for the idealized model of criminal activity introduced in \cite{sbb}.

Let us now comment on previous works. As far as we know, there are three mathematical works related to the crime model (\ref{sysdyn}). Short, Bertozzi and Brantingham \cite{sbb} proposed this model based on mean field considerations. They have also performed a weakly nonlinear analysis on (\ref{sysoriginal}) about the constant solution
$$ (A, \rho)= \left(\gamma +A_0, \frac{\gamma}{\gamma+A_0} \right) $$
assuming that both $A_0 (x)$ and $\gamma (x)$ are homogeneous.
 Rodriguez and Bertozzi have further shown local existence and uniqueness of solutions
 \cite{rb1}.
 In  \cite{ccm}, Cantrell, Cosner and Manasevich have given a rigorous proof of the bifurcations from this constant steady state. On the other hand, in the isotropic case, Kolokolnikov, Ward and Wei  \cite{kww1} have studied existence and stability of multiple symmetric and asymmetric spikes  for (\ref{sysdyn}) using formal matched asymptotics. They derived qualitative results on competition instabilities and Hopf bifurcation and gave some extensions to two-space dimensions.

 The present paper provides rigorous justification for many of the results in \cite{kww1} and also derives some extensions.
 In particular, we establish here the following three new results: first, we reduce the quasilinear chemotaxis problems to a Schnakenberg type reaction-diffusion system and prove the existence of symmetric $k$ spikes. Second,  this paper gives the first rigorous proof of the existence of asymmetric spikes in the isotropic case. Third, we study the pinning effect in an inhomogeneous setting $A_0 (x)$ and $\gamma (x)$. The stability of these spikes is an interesting issue which should be addressed in the future.

We should mention that another model of criminality has been proposed and analyzed by Berestycki and Nadal \cite{bn}. In a forthcoming paper \cite{bw}, we shall study the existence and stability of hotspots (spikes) in this system as well. It is quite interesting to observe that both models admit hotspot (spike) solutions.

The structure of this paper is as follows. We formally construct a one-spike solution in Section 2 in which we state our main results. In Section 3 we show how to approximate the crime hotspot model by the Schnakenberg model. Section 4  is devoted to  the  computation of  the amplitudes and positions of the spikes to leading order. Nondegeneracy conditions are derived in Section 5.
These are required for the existence proof, given in Sections 6--8. In Section 6 we introduce and study the approximate solutions. In Section 7 we apply Liapunov-Schmidt reduction to this problem. Lastly, we solve the reduced problem in Section 8 and conclude the existence proof.  In Section 9 we extend the proof of single spike solution to the case when both $A_0 (x)$ and $\gamma (x)$ are allowed to be inhomogeneous. Finally, in Section 10 we discuss our results and their significance  and mention possible future work and open problems.

\section{Steady state: Formal argument for leading order and main results}
\setcounter{equation}{0}

 Before stating the main results, we first construct a time-independent spike on the interval $[-L,L]$ located
at some point $x_0$. The construction here is carried out using classical matched asymptotic expansions.

In the inner region, we assume that $v$ is a constant $v_0$ in leading order:
\[
v(x)\sim v_{0},\ \ \ |x-x_0|\ll 1.
\]
Then, if $0<\ep\ll 1$, the equation for $A$ becomes
\[ \varepsilon^2 A^{''}-A+ v_0 A^3 + A_0 (x)=0.\]
Rescaling
\[
A(x)=v_{0}^{-1/2} \hat{A}(y),\ \ \ \ y=\frac{x-x_0}{\varepsilon},
\]
we get
\[ \hat{A}_{yy}-\hat{A}+ \hat{A}^3 + A_0 (x_0+\epsilon y) v_0^{1/2}=0.\]
We assume that $v_0 \to 0$ as $\ep\to 0$.
Then, at leading order, $\hat{A} (y) \sim w(y)$, where $w$ is the unique (even) solution of the following ODE
\[
w_{yy}-w+w^{3}=0
\]
so that%
\[
w(y)=\sqrt{2}\operatorname*{sech}\left(  y\right)  .
\]

In the outer region, we assume that%
\[
vA^{3}\ll1,\ \ \ \frac{x}{\varepsilon}\gg1
\]
so that%
\[
A\sim A_{0} (x).%
\]
We also assume that $D=\frac{\hat{D}}{\ep^2}$, where $\hat{D}$ is a positive constant,
and we estimate
\[
\int_{-L}^{L}vA^{3}\,dx\sim v_{0}^{-1/2}\varepsilon\int_{-\infty}^{\infty}w^{3}dy.
\]
Integrating the second equation in (\ref{sysdyn}), we then have
\begin{gather}
\nonumber
v_{0}^{-1/2}\varepsilon\int_{-\infty}^{\infty}w^{3}dy\sim\int_{-L}^L \gamma (x) dx,\\
v_{0}\sim\frac{\left(  \int_{-\infty}^{\infty}w^{3}dy\right)  ^{2}}%
{ \left(\int_{-L}^L \gamma (x) dx \right)^{2}}\varepsilon^{2}.\label{sys0}%
\end{gather}
We remark that $\int w^{3}\,dy=\int w\,dy=\sqrt{2}\pi$ so that%
\[
v_{0}\sim\frac{2 \pi^{2}}{\left(  \int_{-L}^L \gamma (x) \right)  ^{2}}
\varepsilon^{2}.%
\]
In particular, we obtain%
\begin{equation}
A(x)\sim\left\{
\begin{array}
[c]{c}%
A_0(x)+
\dfrac{\sqrt{2} \int_{-L}^L \gamma (x) dx }{\varepsilon\pi}\,w\left(\frac{x-x_0}{\varepsilon}
\right),\ \ \ x=O\left(  \varepsilon\right),  \\[3mm]
A_{0}(x),\ \ \ \ x\gg O(\varepsilon).
\end{array}
\right.  \label{Aunif}%
\end{equation}
Now we state our main theorems on the existence of multi-spike steady states for system (\ref{sysdyn}). We discuss two cases.

 In the case of isotropic coefficients $A_0 (x) \equiv$ Constant, $\gamma (x)\equiv$ Constant, we will consider two types of solutions:

(i) Multiple spikes of arbitrary number having the same amplitude (symmetric spikes).

(ii) Multiple spikes having different amplitude for the case of one large and one small spike (asymmetric spikes).

In the case of anisotropic coefficients $A_0 (x)$ and $\gamma (x)$, we will consider the existence of single spike solution.

\medskip

Our first result concerns the existence of multiple spikes of arbitrary number having the same amplitude (symmetric spikes).

\begin{theorem}
\label{existencesym}
Assume that $D=\frac{\hat{D}}{\ep^2}$ for some fixed $\hat{D}>0$ and
\begin{equation}
\label{aniso}
A_0(x)\equiv A_0, \gamma (x) \equiv \bar{A}-A_0 \ \ \mbox{where}\  \bar{A}>A_0.
\end{equation}
 Then, provided $\ep>0$ is small enough, problem
(\ref{sysdyn})  has a $K$-spike steady state $(A_\ep,v_\ep)$ which satisfies the following properties:
\begin{equation}
\label{aep}
A_{\ep}(x)= A_0+\frac{1}{\ep}\sum_{j=1}^K
\frac{1}{\sqrt{v_j^\ep}}
w \left(\frac{x-t_j^\ep}{\ep}\right)+O\left(\ep\log\frac{1}{\ep}\right),
\end{equation}
\begin{equation}
\label{hep}
v_\ep(t_i^\ep)=\ep^2
v_i^\ep
,\quad
i=1,\ldots,K,
\end{equation}
where
\begin{equation}
\label{tep0}
t_i^\ep\to t_i^0,\quad
i=1,\ldots,K
\end{equation}
with
\begin{equation}
\label{limpos}
t_i^0=\frac{2i-1-K}{K}\,L,\,i=1,\ldots,K
\end{equation}
and
\begin{equation}
\label{vep0}
v_i^\ep= v_i^0\left(1+O\left(\ep\log\frac{1}{\ep}
\right)\right)
,\quad
i=1,\ldots,K
\end{equation}
with
\begin{equation}
\label{limamp}
v_i^0=\frac{\pi^2 K^2}{2 (\bar{A}-A_0)^2 L^2},\quad
i=1,\ldots,K.
\end{equation}
\end{theorem}

\begin{remark} {\rm
{\rm 
Note that in (\ref{aep}) a two-term expansion of the solution $A_\ep$ is given, where for each spike the term
$\frac{1}{\sqrt{v_j^\ep}}
w \left(\frac{x-t_j^\ep}{\ep}\right)$
of order $O(\frac{1}{\ep})$ is the leading term in the inner solution and the term
$A_0$ of order $O(1)$ is the leading term in the outer solution.
By using the operator $T[\hat{A}]$ defined in (3.12) this two-term expansion
 carries over to $\hat{v}$ as well.
The same remark applies to (\ref{aep1}) and (\ref{aep200}).
The two-term expansion agrees with that in \cite{kww1}.
} 
} \end{remark}

The next result is about asymmetric two-spikes.

\begin{theorem}
\label{existenceas}
Under  the same assumption as in Theorem \ref{existencesym}, with $D=\frac{\hat{D}}{\ep^2}$ for some fixed $\hat{D}>0$ and suppose moreover that
\begin{equation}
\frac{2\sqrt{\pi}(\hat{D} A_0^2)^{1/4}}{(\bar{A}-A_0)^{3/4}L}
\leq 1,
\label{condc}
\end{equation}
and
\begin{equation}
\frac{2\sqrt{\pi}(\hat{D} A_0^2)^{1/4}}{(\bar{A}-A_0)^{3/4}L}
\not= \frac{2}{\sqrt{5}}.
\label{condbc}
\end{equation}
Then, for $\ep>0$ small enough, problem
(\ref{sysdyn})  has an asymmetric $2$-spike steady state $(A_\ep,v_\ep)$ which satisfies the following properties:
\begin{equation}
\label{aep1}
A_{\ep}(x)=A_0+\frac{1}{\ep}
\left(
\sum_{j=1}^2
\frac{1}{\sqrt{v_i^\ep}}
w \left(\frac{x-t_i^\ep}{\ep}\right)
+O\left(\ep\log\frac{1}{\ep}\right)
\right),
\end{equation}
\begin{equation}
\label{hep1}
v_\ep(t_i^\ep)=\ep^2
v_i^\ep,\quad
i=1,\ldots,K,
\end{equation}
where $t_i^\ep$ and $v_i^\ep$ satisfy (\ref{tep0}) and (\ref{vep0}), respectively.
The limiting amplitudes $v_i^0$ and positions $t_i^0$ are given as solutions of (\ref{amp1}) and (\ref{amp6}).
\end{theorem}

Condition (\ref{condbc}) is a kind of nondegeneracy condition.
Note that in the case of asymmetric spikes we explicitly characterize the points of non-degeneracy.

The last theorem is about the existence of single spike solution in the anisotropic case

\begin{theorem}
\label{existenceani}
Assume that
$\ep>0$ is small enough and
$D=\frac{\hat{D}}{\ep^2}$ for some fixed $\hat{D}>0$. Then, problem
(\ref{sysdyn})  has a single spike steady state $(A_\ep,v_\ep)$ which satisfies the following properties:
\begin{equation}
\label{aep200}
A_{\ep}(x)= A_0 (x)+\frac{1}{\ep}
\frac{1}{\sqrt{v_0^\ep}}
w \left(\frac{x-t_0^\ep}{\ep}\right)+O\left(\ep\log\frac{1}{\ep}\right),
\end{equation}
\begin{equation}
\label{hep200}
v_\ep(t_0^\ep)=\ep^2
v_0^\ep,
\end{equation}
where
\begin{equation}
\label{tep000}
t_0^\ep\to t_0, \ \int_{-L}^{t_0} \gamma (x) dx = \int_{t_0}^L \gamma (x) dx
\end{equation}
and
\begin{equation}
\label{vep000}
v_0^\ep= \frac{ 2\pi^2}{ (\int_{-L}^L \gamma (x) dx)^2} \left(1+  O\left(\ep\log\frac{1}{\ep}\right)\right).
\end{equation}
\end{theorem}
We notice that in the anisotropic case, the single spike location is only determined by the function $\int_{-L}^{x} \gamma (t) dt$ and $A_0 (x)$ has no effect at all. Note also that the location $t_0$ is uniquely determined by the condition
\begin{equation}
\int_{-L}^{t_0} \gamma (x) dx= \frac{1}{2} \int_{-L}^L \gamma (x) dx.
\end{equation}
With more computations, it is possible  to construct multiple asymmetric spikes in the isotropic case, and also multiple spikes in the anisotropic case. Since the statements and computations are complicated, we will not present them here. We refer to \cite{ww-pre} for some results in this direction.

\section{Scaling and approximation by the Schnakenberg model}
\setcounter{equation}{0}

We will use the following notation
for the domain and the rescaled domain, respectively:
\begin{equation}
 \Omega= (-L, L),\quad \Om_{\ep}=\left(-\frac{L}{\ep},\frac{L}{\ep}\right).
 \end{equation}

This section is devoted to the reduction of the system (\ref{sysdyn}) to a particular Schnakenberg type reaction diffusion equation in which no chemotaxis appears.

Based on the computations in Section 2,  we rescale the solution and the second diffusion coefficient as follows:
\[ A=A_0 (x)+\frac{1}{\epsilon} \hat{A},\quad v=\epsilon^2 \hat{v},\ \ \  D=\frac{\hat{D}}{\epsilon^2}.\]

Then the steady-state problem becomes
\begin{align}
\nonumber
0  & =\varepsilon^{2}\hat{A}_{xx}-\hat{A}+\hat{v} (\epsilon A_0 +\hat{A})^{3} +\varepsilon^3 A_0^{''},\ x\in\Om, \\
0 & =\hat{D} \left(  \left(A_0 (x)+\frac{1}{\epsilon} \hat{A}\right)^2 \hat{v}_{x}\right)  _{x}-\frac{1}{\epsilon} \hat{v} (\epsilon A_0(x) +\hat{A})^{3}+\gamma (x),\ x\in\Om.
\label{sys1}
\end{align}

We will consider the case when $\ep\ll 1$ and $\hat{D}$ is constant, with Neumann boundary conditions.

A key observation of this paper is that the solutions of
problem (\ref{sys1}) are very close to the solutions of the Schnakenberg model
\begin{align}\nonumber
0  & =\varepsilon^{2}\hat{A}_{xx}-\hat{A}+\hat{v} (\epsilon A_0 +\hat{A})^{3} +\varepsilon^3 A_0^{''},\ x\in\Om,\\
0 & =\hat{D} \left(  A_0^2 \hat{v}_{x}\right)  _{x}-\frac{1}{\epsilon} \hat{v} (\epsilon A_0 +\hat{A})^{3}+\gamma (x),\ x\in\Om,
\label{schnak1}
\end{align}
with Neumann boundary conditions.

To see this, we first consider the following linear problem:
\begin{equation}
\label{vequ}
\left\{
\begin{array}{l}
\hat{D} (a(x) v_x)_x=f(x),\quad -L<x<L,
\\[3mm]
v_x(-L)=v_x(L)=0,
\end{array}
\right.
\end{equation}
where $a\in C^1(-L,L)$, $a(x)\geq c>0$ for all $x\in(-L,L)$ and $f\in L^1(-L,L)$.
We compute
\[ a(x) v_x(x)=\int_{-L}^x \frac{1}{ \hat{D}} f(t)\,dt.\]
So
\begin{equation}
\label{v1}
v_x (x)=\frac{1}{ a(x)} \int_{-L}^x \frac{1}{ \hat{D}} f(t)\,dt
\end{equation}
and hence
\begin{equation}
\label{ss2}
v(x)- v(-L)= \int_{-L}^x \frac{1}{ \hat{D} a(s)} \int_{-L}^s f(t) dt\,ds
\end{equation}
which can be rewritten as
\begin{equation}
\label{ss3}
v(x)- v(-L)=\frac{1}{\hat{D}} \int_{-L}^x K_a (x, s) f(s)\, ds,
\end{equation}
where
\[ K_a (x, s)=\int_s^x \frac{1}{a (t)} \,dt.\]
\begin{remark} {\rm
We note that the kernel $K_a (x, s)$
is an even (odd) function if $a(x)$ is an odd (even) function. More precisely, if $a(x)=\pm a(-x)$,  then
\begin{equation}
\label{ksym}
K_a (-x, -s)=\mp K_a (x,s).
\end{equation}
} \end{remark}

\begin{remark} {\rm
We note that $v$ is an even (odd) function if $f$ is even and $a$ is even (odd). More precisely, using
 \[
 \int_{-L}^0 f(t)\,dt=\frac{1}{2}\int_{-L}^L f(t)\,dt=0 \quad \mbox{if $f$ is even},
 \]
we compute
\[
v_x(x)=\frac{1}{a(x)}\int_0^x \frac{1}{\hat{D}}f(t)\,dt.
\]
Integration yields
 \[
v(x)-v(0)=\int_0^x \frac{1}{\hat{D}a(s)}\int_0^s f(t)\,dt\,ds
\]\[
=
\frac{1}{\hat{D}}
\int_{0}^x K_a(x,s)f(s)\,ds
\]
and
\[
v(-x)-v(0)=\frac{1}{\hat{D}}
\int_{0}^{-x} K_a(-x,s)f(s)\,ds
\]
\[
=-\frac{1}{\hat{D}}
\int_{0}^{x} K_a(-x,-s)f(-s)\,ds
\]
\[
=\pm\frac{1}{\hat{D}}
\int_{0}^{x} K_a(x,s)f(s)\,ds
\]
\[
=\pm(v(x)-v(0))
\]
if $a$ is an even (odd) function
using (\ref{ksym}).
Similarly, if $f$ is odd and $a$ is odd (even), then $v$ is an even (odd) function.
} \end{remark}

Integrating (\ref{vequ}), we derive the necessary condition
\begin{equation}
\label{intconstr}
\int_{-L}^L f(x)\,dx=0.
\end{equation}
Note that on the other hand $v$ defined by (\ref{vequ}) satisfies the boundary conditions
$v_x(-L)=v_x(L)=0$  provided that (\ref{intconstr}) holds. This follows from (\ref{v1}).

Let us now consider $a(x)= \left(A_0+ \frac{\gamma }{\epsilon} w(\frac{x}{\epsilon})\right)^2$, where $ w>0$ and $ w(y)\sim e^{-|y|}$ as $|y|\to\infty$. Then we claim that
\begin{equation}
\label{ss4}
K_{a} (x,s)= K_{A_0^2}  (x,s)+ O(\epsilon |s-x|)+ O\left( \left|[s,x] \cap \left(-2\epsilon \log\frac{1}{\epsilon}, 2\epsilon \log\frac{1}{\epsilon}\right)\right|\right).
\end{equation}
Note that (\ref{ss4}) is an $L^\infty$ estimate for $K_a(x,s)$.

In fact, we have
\[ \int_{s}^x \frac{1}{(A_0+\frac{1}{\epsilon} w)^2}\,dt = \int_s^x \frac{1}{A_0^2} dx+ \int_{s}^x \left[\frac{1}{(A_0+\frac{1}{\epsilon} w)^2}-\frac{1}{A_0^2} \right]\,dt,\]
where
\[ \int_{s}^x \left[\frac{1}{A_0^2}-\frac{1}{(A_0+\frac{1}{\epsilon} w)^2}\right]\,dt= \epsilon \int_{\frac{s}{\epsilon}}^{\frac{x}{\epsilon}} \frac{ 2\epsilon A_0 w+w^2}{ (\epsilon A_0+w)^2}\,dy\]
\[= \epsilon \int_{ [s/\epsilon, x/\epsilon] \cap \left\{ |y|>2 \log \frac{1}{\epsilon}\right\} }
\ldots \,dy + \epsilon
\int_{ [s/\epsilon, x/\epsilon] \cap \left\{ |y|<2 \log \frac{1}{\epsilon}\right\} }\ldots \,dy. \]
The first term is $O(\epsilon |x-s|)$ since $w=O(\ep^2)$ and so $ \frac{2\epsilon A_0 w+w^2}{ (\epsilon A_0+w)^2 }=O(\epsilon)$. For the second term, observing that $ \frac{2\epsilon A_0 w+w^2}{ (\epsilon A_0+w)^2 }=O(1)$ we derive (\ref{ss4}). All these estimates are in the $L^\infty$ norm.

Thus, $v$ satisfies:
\begin{equation}
\label{ss5}
v(x)- v(-L)=\frac{1}{\hat{D}} \int_{-L}^x K_{A_0^2} (x, s) f(s) ds
+O\left(\epsilon \int_{-L}^x \left(|x-s|+\log\frac{1}{\ep}\right)\, |f(s)|\, ds\right).
\end{equation}

\begin{remark} {\rm
The estimates (\ref{ss4}) and (\ref{ss5}) also hold if
$$ a (x)= \left(A_0 +\frac{\gamma}{\epsilon} \left(w\left(\frac{x-x_0}{\epsilon}\right) + \phi\right) \right)^2,$$
where $ \phi (x)$ satisfies $ |\phi (x)|\leq C \epsilon \max ( e^{- \frac{|x-x_0|}{2\epsilon}}, \sqrt{\epsilon})$.  This is the class of functions that we will work with. This is also the motivation for our choice of the norm $\| \cdot \|_{*}$ (defined in (\ref{normdef})).
} \end{remark}

Therefore, we can approximate steady states for the crime hotspot model by the Schnakenberg model as follows: Given $\hat{A}>0$, let $\hat{v}=T[\hat{A}]$ be the unique solution of the following linear problem:
\begin{equation}
\left\{
\begin{array}{l}
\hat{D} \left(  (A_0 +\frac{1}{\epsilon} \hat{A})^2 \hat{v}_{x}\right)  _{x}-\frac{1}{\epsilon} \hat{v} (\epsilon A_0 +\hat{A})^{3}+\gamma (x)=0,
\quad -L<x<L,
\\[3mm]
\hat{v}_{x}(-L)=\hat{v}_{x}(L)=0.
\end{array}
\right.
\label{ta}
\end{equation}
Then, by the maximum principle, the solution $T[\hat{A}]$ is positive.

By the previous computations and remarks,  if $ \hat{A}=  w+\phi$ with $ |\phi |\leq C \epsilon \max ( e^{-\frac{|x-x_0|}{2\epsilon}}, \sqrt{\epsilon})$, it follows that
\begin{equation}
\label{approx}
T[\hat{A}]= v^{0} +O\left(\epsilon \log \frac{1}{\epsilon}\right) \quad\mbox{ in }H^2(-L,L),
\end{equation}
where $v^0$ satisfies
\begin{equation}
\left\{
\begin{array}{l}
\hat{D} \left(  A_0^2 v^0_{x}\right)  _{x}-\frac{1}{\epsilon} v^0 (\epsilon A_0 +\hat{A})^{3}+\gamma (x)=0,
\quad
-L<x<L,
\\[3mm]
v^0_{x}(-L)=v^0_{x}(L)=0.
\end{array}
\right.
\end{equation}


We adapt an approach based on Liapunov-Schmidt reduction which has been applied to the semilinear Schnakenberg model in \cite{iww2} and extend it to the quasilinear crime hotspot model.  This method has also been used to study spikes for the one-dimensional Gierer-Meinhardt system in \cite{ww, ww-pre} as well as  two-dimensional Schnakenberg model in \cite{ww13}. We refer to the survey paper \cite{wei-survey} and the book \cite{ww-book} for references. Multiple asymmetric spikes for the one-dimensional Schnakenberg model have been considered using matched asymptotics in \cite{iww}. Existence and stability of localized patterns for the crime hotspot model have been studied by matched asymptotics in \cite{kww1} and results on competition instabilities and Hopf bifurcation have been shown including some extensions to two space dimensions.

We  remark that another approach  for studying multiple spikes in one-dimensional reaction-diffusion system is the geometric singular perturbation theory in dynamical systems. For results and methods in this direction we refer to \cite{dgk, dkp} and the references therein.

\section{Computation of the amplitudes and positions of the spikes}
\setcounter{equation}{0}

In this section, we study (\ref{sys1}) in the isotropic case (\ref{aniso}).
 In particular, we
compute the amplitudes and positions to leading order. We consider symmetric multi-spike solutions with any number of spikes and asymmetric multi-spike solutions with one small and large spike.

We first write down the system for the amplitudes in case of a general number $K$ of spikes, where we have either $K$ spikes of the same amplitude or $k_1$ small and $k_2$ large spikes with $k_1+k_2=K$. We will first solve this system in the case of symmetric spikes. Then we will choose $k_1=k_2=1$ and solve this system in this special case of asymmetric spikes.

Integrating the right hand side of the second equation in (\ref{sys1}), we compute for $v_j=\lim_{\ep\to 0} \hat{v}_\ep(t_j^\ep)$:
\begin{equation}
\label{amp1}
\sum_{j=1}^K
\frac{\sqrt{2}\pi}{\sqrt{v_j}}=
(\bar{A}-A_0)2L.
\end{equation}
Solving the second equation in (\ref{sys1}), using (\ref{ss2})
in combination with the approximation (\ref{approx}), we get
\[
\hat{v}_\ep(x)- \hat{v}_\ep(-L)= \int_{-L}^x \frac{1}{ \hat{D} (A_0+\frac{1}{\ep}\hat{A}_\ep)^2} \int_{-L}^s
\left(
\frac{1}{\ep}\hat{v}_\ep(\ep A_0+\hat{A}_\ep)^3-\bar{A}+A_0
\right)
 dt\,ds
\]\[
=\int_{-L}^x \frac{1}{\hat{D} A_0^2} \int_{-L}^s
\left(
\frac{1}{\ep}\hat{v}_\ep(\ep A_0+\hat{A}_\ep)^3-\bar{A}+A_0
\right)
 dt\,ds
 +O\left(\ep\log\frac{1}{\ep}\right)
\]\[
\int_{-L}^x \frac{1}{\hat{D} A_0^2}
\left[
\sum_{j=1}^{i-1} \frac{\sqrt{2}\pi}{\sqrt{\hat{v}_\ep(t_j^\ep)}}
H(s-t_j^\ep)-
\left(\sum_{j=1}^K \frac{\sqrt{2}\pi}{\sqrt{\hat{v}_\ep(t_j^\ep)}}
\right)
\frac{s+L}{2L}
\right]
\,ds
+O\left(\ep\log\frac{1}{\ep}\right)
\]\[
=
\frac{1}{\hat{D} A_0^2}
\left[
\sum_{j=1}^{i-1} \frac{\sqrt{2}\pi}{\sqrt{\hat{v}_\ep(t_j^\ep)}}
(x-t_j^\ep)-
\left(\sum_{j=1}^K \frac{\sqrt{2}\pi}{\sqrt{\hat{v}_\ep(t_j^\ep)}}\right)
\frac{(x+L)^2}{4L}
\right]
+O\left(\ep\log\frac{1}{\ep}\right),
\]
where
\[
H(x)=
\left\{
\begin{array}{l}
1 \quad \mbox{ if } x\geq 0,
\\[3mm]
0 \quad \mbox{ if }x<0.
\end{array}
\right.
\]
Taking the limit $\ep\to 0$ and setting $x=t_i=\lim_{\ep\to 0}t_i^\ep$, we derive
\[
v_i=
\frac{1}{\hat{D} A_0^2}
\left[
\sum_{j=1}^{i-1} \frac{\sqrt{2}\pi}{\sqrt{v_j}}(t_i-t_j)-
\left(\sum_{j=1}^K \frac{\sqrt{2}\pi}{\sqrt{v_j}}\right)
\frac{(t_i+L)^2}{4L}
\right]+C_1
\]
\begin{equation}
\label{amp2}
=\frac{1}{\hat{D}A_0^2}\left[
\sum_{j=1}^K
\frac{\sqrt{2}\pi}{\sqrt{v_j}}\frac{1}{2}
|t_i-t_j|
-\left(\sum_{j=1}^K
\frac{\sqrt{2}\pi}{\sqrt{v_j}}\right)
\frac{1}{4L}t_i^2
\right]+C_2
\end{equation}
for some real constants $C_1,\,C_2$ independent of $i$, where the last identity in (\ref{amp2}) uses (\ref{amp3}) which we now explain.


We use an assumption on the position of spikes that can be stated as follows:
\begin{equation}
\label{fi}
F_i(t_1^0,t_2^0,\ldots,t_K^0):=
\frac{1}{2}
\frac{\sqrt{2}\pi}{\sqrt{v_i}}
+\sum_{j=1}^{i-1} \frac{\sqrt{2}\pi}{\sqrt{v_j}}-
\left(
\sum_{j=1}^K \frac{\sqrt{2}\pi}{\sqrt{v_j}}
\right)
\frac{t_i+L}{2L}
=0,\quad i=1,\ldots,K.
\end{equation}
Note that (\ref{fi}) will be derived later on in
Section 8 below (see equation (\ref{formulaf})).
We re-write (\ref{fi}) and compute
\[
\frac{1}{2\sqrt{v_i}}
+\sum_{j=1}^{i-1}
\frac{1}{\sqrt{v_j}}
-\left(\sum_{j=1}^{K}
\frac{1}{\sqrt{v_j}}\right)
\frac{t_i+L}{2L}
\]
\begin{equation}
\label{amp3}
=\sum_{j=1}^{i-1}
\frac{1}{2\sqrt{v_j}}
-
\sum_{j=i+1}^{K}
\frac{1}{2\sqrt{v_j}}
-
\left(
\sum_{j=1}^{K}
\frac{1}{\sqrt{v_j}}
\right)
\frac{t_i}{2L}=0.
\end{equation}


From (\ref{amp2}) and (\ref{amp3}), we derive
\begin{equation}
\label{amp5}
v_{i+1}-v_i=\frac{\pi}{\sqrt{2} \hat{D}A_0^2}
(t_{i+1}-t_i)\frac{1}{2}
\left(\frac{1}{\sqrt{v_i}}-
\frac{1}{\sqrt{v_{i+1}}}
\right).
\end{equation}
In the next two subsection we now solve these equations for the amplitudes of the spikes in the cases of both symmetric and asymmetric spikes.

\subsection{Symmetric spikes}

We first consider the case of symmetric spikes, where
$v_i=v$ is independent of $i=1,\ldots,K$,
and compute the amplitude $v$ and the positions $t_i$.

From (\ref{amp1}), we get
\[
v=\frac{\pi^2 K^2}{2(\bar{A}-A_0)^2  L^2}.
\]

From (\ref{amp5}), the positions are
$t_i=\left(-1+\frac{2i-1}{K}\right)L,\quad i=1,\ldots,K.$

The proof of the existence of multiple symmetric spikes follows from the construction of a single spike in the interval
$\left(-\frac{L}{K},\frac{L}{K}\right)$. The proof of the existence of a single spike uses the implicit function theorem in the space of even functions, for which the Liapunov-Schmidt reduction method is not needed.
This proof can easily be obtained by specializing  the proof given for multiple asymmetric  spikes given below. In this case, the proof can thus be simplified.
Therefore we omit the details.

\subsection{Asymmetric spikes}

Combining (\ref{amp3}) and (\ref{amp5}), we get
\[
v_iv_{i+1}=\frac{\pi}{\sqrt{2}\hat{D}A_0^2}
\frac{L}{2} \frac{1}{\sum_{j=1}^K \frac{1}{\sqrt{v_j}}}.
\]
This implies that there are only two different amplitudes which we denote by
$v_s\leq v_l$ appearing $k_1$ and $k_2$ times, respectively. Hence we get
\begin{equation}
\label{amp6}
v_sv_l=\frac{\pi}{\sqrt{2}\hat{D}A_0^2}
\frac{L}{2} \frac{1}{ \frac{k_1}{\sqrt{v_s}}
+
\frac{k_2}{\sqrt{v_l}}
}.
\end{equation}

Multiplying (\ref{amp3}) by  $\frac{t_i}{2}$ and subtracting (\ref{amp2}) from the result
we get
\[v_i-C=\frac{\pi}{\sqrt{2}\hat{D}A_0^2}
\sum_{j=1}^k\frac{t_j}{\sqrt{v_j}}\mbox{sgn}(t_j-t_i)
+\frac{1}{\hat{D}A_0^2}
\left(\sum_{j=1}^k\frac{\pi}{\sqrt{2}}\frac{1}{\sqrt{v_j}}\right)
\frac{t_i^2}{2L},
\]
where
\[
\mbox{sgn}(\alpha)=\left\{
\begin{array}{ll}
+1 & \mbox{ if }\alpha>0,
\\[2mm]
0 & \mbox{ if }\alpha=0,
\\[2mm]
-1 & \mbox{ if }\alpha<0.
\end{array}
\right.
\]

Next we determine $v_s,\,v_l$ from
(\ref{amp1}) and (\ref{amp6}). Substituting (\ref{amp1}) into (\ref{amp6}), we get
\[
v_l=\frac{1}{v_s}\frac{\pi^2}{\hat{D}A_0^2}\frac{1}{4}
\frac{1}{\bar{A}-A_0}.
\]
Plugging this equation into (\ref{amp1}) gives
\[
C\left(z+\frac{1}{z}\right)=1,
\]
where
\[ C
=\frac{\sqrt{\pi}(\hat{D}A_0^2)^{1/4}}{(\bar{A}-A_0)^{3/4}L}
\sqrt{k_1k_2}, \quad
z=
\frac{\sqrt{2v_s}(\hat{D}A_0^2)^{1/4}(\bar{A}-A_0)^{1/4}}{\sqrt{\pi}}
\sqrt{\frac{k_2}{k_1}}.
\]
To determine a solution, we need to satisfy the necessary condition $2C<1$
which can be summarized as
\[
\frac{2\sqrt{\pi}(\hat{D}A_0^2)^{1/4}}{(\bar{A}-A_0)^{3/4}L}
\sqrt{k_1k_2}<1.
\]
The second necessary condition
is given by
$v_s<v_l$ which is equivalent to $z< \sqrt{\frac{k_2}{k_1}}$.
This implies the following cases:

\medskip

\noindent
{\bf Case (i):} $k_2\leq k_1$.

If
\[C\left(\sqrt{\frac{k_2}{k_1}}+\sqrt{\frac{k_1}{k_2}}\right) <  1 \]
then there exists exactly one solution with $z\leq  \sqrt{\frac{k_2}{k_1}}$.

On the other hand, if
\[C\left(\sqrt{\frac{k_2}{k_1}}+\sqrt{\frac{k_1}{k_2}}\right)> 1 \]
then there exists no solution with $z<\sqrt{\frac{k_2}{k_1}}$.

\medskip
\noindent
{\bf Case (ii):} $k_2> k_1$.

If $2C>1$ there is no solution. If $2C<1$ and \[C\left(\sqrt{\frac{k_2}{k_1}}+\sqrt{\frac{k_1}{k_2}}\right)< 1 \]
then there exists exactly one solution with $z\leq \sqrt{\frac{k_2}{k_1}}$.

If $2C<1$ and \[C\left(\sqrt{\frac{k_2}{k_1}}+\sqrt{\frac{k_1}{k_2}}\right)> 1 \]
then there exist exactly two solutions with $z< \sqrt{\frac{k_2}{k_1}}$.

\medskip

\noindent
{\bf Special case:} $k_1=k_2=1$.

Finally, we consider the special case $k_1=k_2=1$ which belongs to Case (i) in the previous classification and we have the following results:

If $2C< 1$ then there is one solution with $z< 1. $
If $2C> 1$ then there is no solution with $z> 1. $


\section{Existence and nondegeneracy conditions}
\setcounter{equation}{0}

We now describe a general scheme of Liapunov-Schmidt reduction. (We refer to the survey paper \cite{wei-survey} for more details.) Essentially this method divides the problem of solving nonlinear elliptic equations (and systems) into two steps. In the first step, the problem is solved up to multipliers of approximate kernels. In the second step one solves algebraic equations in terms of finding zeroes of the multipliers.



In this section, we  linearize (\ref{sys1}) around the approximate solution and derive the linearized operator as well as its nondegeneracy conditions, i.e. conditions such that the resulting linear operator is uniformly invertible.

\medskip

Linearizing (\ref{sys1}) around the solution, we get:
\begin{align}\nonumber
0  & =\varepsilon^{2}\phi_{xx}-\phi+3\hat{v}(\ep A_0+\hat{A})^2\phi+\psi (\ep A_0+\hat{A})^3\ \ \mbox{ in }\Om, \\[3mm]
\nonumber
0 & =\hat{D}
\left( \left(A_0+\frac{1}{\ep}\hat{A}\right)^2
\psi_{x}\right)_x
+\hat{D}
\left(2\left(A_0+\frac{1}{\ep}\hat{A}\right)
\frac{1}{\ep}\phi
\hat{v}_{x} \right)_x
\\[3mm] &
-\frac{3}{\epsilon} \hat{v}(\ep A_0+\hat{A})^2\phi -\frac{1}{\epsilon} \psi (\ep A_0+\hat{A})^3 \ \ \mbox{ in }\Om
\label{evp}%
\end{align}
with Neumann boundary conditions
\[
\phi_x(-L)=\phi_x(L)=\psi_x(-L)=\psi_x(L)=0.
\]
Note that for the second equation of (\ref{evp}) we have the necessary condition
\[
\int_{-L}^L
\left(-\frac{3}{\epsilon} \hat{v}(\ep A_0+\hat{A})^2\phi -\frac{1}{\epsilon} \psi (\ep A_0+\hat{A})^3
\right)\,dx=0
\]
which follows by integrating the equation  and using the Neumann boundary conditions for
$\psi$ and $v$. In the limit $\ep\to 0$ we get
\[
\sum_{j=1}^K \left[\psi_j
\frac{\sqrt{2}\pi}{v_j^{3/2}}
+3\int_{-\infty}^\infty w^2\phi_j\,dy
\right]=0,
\]
where $\psi_j=\psi(x_j)$.

The second equation of (\ref{evp}) can be solved as follows, using formula (\ref{ss2}) and estimate (\ref{approx}):
\[
\psi(x)- \psi(-L)= \int_{-L}^x \frac{1}{ \hat{D} (A_0+\frac{1}{\ep}\hat{A})^2} \int_{-L}^s
\Bigg[
-\hat{D}
\left(
2\left(A_0+\frac{1}{\ep}\hat{A}\right)
\frac{1}{\ep}\phi
\hat{v}_{x}
\right)_x
\]\[
+\frac{3}{\ep}\hat{v}(\ep A_0+\hat{A})^2\phi
+\frac{1}{\ep}\psi(\ep A_0+\hat{A})^3
\Bigg]
 dt\,ds
\]\[
=\int_{-L}^x \frac{1}{\hat{D} A_0^2} \int_{-L}^s
\left(
\frac{3}{\ep}\hat{v}(\ep A_0+\hat{A})^2\phi
+\frac{1}{\ep}\psi(\ep A_0+\hat{A})^3
\right)
 dt\,ds
 +O\left(\ep\log\frac{1}{\ep}\right).
\]
Note that the contributions from the term
\[
\hat{D}
\left(
2\left(A_0+\frac{1}{\ep}\hat{A}\right)
\frac{1}{\ep}\phi
\hat{v}_{x}
\right)_x
\]
can be estimated by $O\left(\ep\log\frac{1}{\ep}\right)$
since $\phi$ vanishes in the outer expansion.
In the limit $\ep\to 0$, we get
\[
\psi(x)- \psi(-L)=
\int_{-L}^x \frac{1}{\hat{D} A_0^2}
\left[
\sum_{j=1}^{i-1} \psi_j \frac{\sqrt{2}\pi}{v_j^{3/2}}H(s-t_j)
+
3\sum_{j=1}^{i-1}
\int_{-\infty}^\infty
 w^2\phi_j\,dy
H(s-t_j)
\right]
\,ds
\]
\begin{equation}
=\frac{1}{\hat{D} A_0^2}
\left[\sum_{j=1}^{i-1} \psi_j \frac{\sqrt{2}\pi}{v_j^{3/2}}(x-t_j)
+
3\sum_{j=1}^{i-1}
\int_{-\infty}^\infty
 w^2\phi_j\,dy\,
(x-t_j)\right],
\label{limeig}
\end{equation}
where $t_{i-1}< x\leq t_{i}$.

From now on, we consider the case of two spikes
having different amplitudes (asymmetric spikes).
Using the notation
\[\Phi=
\left(\begin{array}{l}\phi_1 \\[2mm] \phi_2 
\end{array}
\right),\quad
\Psi=
\left(\begin{array}{l}\psi_1 \\[2mm] \psi_2 
\end{array}
\right),\quad
\omega=
\left(\begin{array}{l}\omega_1 \\[2mm] \psi_2 
\end{array}
\right)
=
\left(\begin{array}{l}-\frac{3}{\hat{D}A_0^2}\int
w^2\phi_1\,dy \\[2mm] -\frac{3}{\hat{D}A_0^2}\int
w^2\phi_2\,dy
\end{array}
\right),
\]
we can rewrite (\ref{limeig}) for $x=t_i$ as follows:
\[
({\cal B}+{\cal C})\Psi=\boldmath{\omega},
\]
where
\[
{\cal C}=\frac{\sqrt{2}\pi}{\hat{D}A_0^2}
\left(
\begin{array}{cc}
\displaystyle
\frac{1}{v_s^{3/2}} & 0 \\[3mm]
\displaystyle
0 &  \displaystyle\frac{1}{v_l^{3/2}}
\end{array}
\right),
\quad
{\cal B}=\frac{1}{d_2}
\left(
\begin{array}{rr} 1 & -1\\
-1 & 1
\end{array}
\right), \quad d_2=t_2-t_1.
\]
Using ${\cal E}={\cal C}({\cal B}+{\cal C})^{-1}$, we get the following system of nonlocal eigenvalue problems (NLEPs)
\begin{equation}
L\Phi:=\Phi_{yy}-\Phi+3w^2\Phi-3w^3\frac{\int w^2{\cal E}\Phi\,dy}{\int w^3\,dy}.
\label{vecnlep}
\end{equation}
Diagonalizing the matrix $\cal E$, we know from \cite{wei99,wz} that (\ref{vecnlep}) has a nontrivial solution iff ${\cal E}$ has eigenvalue $\lambda_e=\frac{2}{3}$.

Thus it remains to compute the matrix ${\cal E}$ and its eigenvalues.

We get
\[
{\cal E}^{-1}=({\cal B}+{\cal C}){\cal C}^{-1}=
{\cal B}{\cal C}^{-1}+I
\]
\[
=\frac{\hat{D}A_0^2}{\sqrt{2}\pi d_2}
\left(
\begin{array}{rr}
v_s^{3/2} & -v_l^{3/2} \\
-v_s^{3/2} & v_l^{3/2}
\end{array}
\right)
+I.
\]
Then ${\cal E}^{-1}$ has the eigenvector $v_{m,1}=\frac{1}{\sqrt{2}}(1,-1)^T$ with eigenvalue
$e_{m,1}=\frac{\hat{D}A_0^2}{\sqrt{2}\pi d_2}
(v_s^{3/2}+v_l^{3/2})+1$
and the eigenvector
$v_{m,2}=\frac{1}{\sqrt{v_s^3+v_l^3}}(v_l^{3/2},v_s^{3/2})^T$ with eigenvalue
$e_{m,2}=1\not=\frac{3}{2}$.

For nondegeneracy, the condition
$e_{m,1}\not=\frac{3}{2}$ has to be satisfied,
which is equivalent to
\[
\frac{\hat{D}A_0^2}{\sqrt{2}\pi d_2}
(v_s^{3/2}+v_l^{3/2})\not= \frac{1}{2}.
\]
Using the formulas for $d_2,\,v_s,\,v_l$, we compute
\[
\frac{\hat{D}A_0^2}{\sqrt{2}\pi d_2}(v_s^{3/2}+v_l^{3/2})
\]\[
=\frac{\hat{D}A_0^2}{\sqrt{2}\pi L}\left((\sqrt{v_s}+\sqrt{v_l})^3-3(\sqrt{v_s}v_l+v_s\sqrt{v_l})\right)
\]\[
=\frac{\hat{D}A_0^2}{\sqrt{2}\pi L}\left(\frac{(\bar{A}-A_0)^{3/2} L^3}{2\sqrt{2}(\hat{D}A_0^2)^{3/2}}
-3 \frac{\pi L}{2\sqrt{2}\hat{D}A_0^2}\right)
\]\[
=\frac{(\bar{A}-A_0)^{3/2} L^2}{4\pi\sqrt{\hat{D}A_0^2}}-\frac{3}{4}\not=\frac{1}{2}.
\]
This implies the condition
\[
\frac{2\sqrt{\pi}(\hat{D}A_0^2)^{1/4}}{(\bar{A}-A_0)^{3/4}L}\not= \frac{2}{\sqrt{5}}.
\]

We have to exclude this point from our existence result Theorem \ref{existenceas}. This is why we impose  the condition (\ref{condbc}) in Theorem \ref{existenceas}, which amounts to a nondegeneracy condition. If this condition is violated, we expect small eigenvalues to occur and it is an open question to know whether there will be spikes in this case.

\section{Approximate solutions}
\setcounter{equation}{0}

For simplicity, we set $L=1$. In this section and the following we consider the case of general $K=1,2,\ldots$ since it does not cause any extra difficulty here, even in the case of asymmetric spikes.
Let $ -1<t_1^0 < \cdots < t_j^0 < \cdots t_K^0 <1$ be $K$ points
and let
$v_j>0$ be $K$ amplitudes
satisfying the assumptions
(\ref{amp1}), (\ref{amp2}) and (\ref{amp3}). Let
\begin{equation}
{\bf t}^0=(t_{1}^0,\ldots,t_{K}^0).
\end{equation}

\medskip

We first construct an approximate solution to (\ref{sys1}) which
concentrates near these prescribed $K$ points.
Then we will rigorously construct an exact solution which is given by a small perturbation of this approximate solution.

Let $-1<t_1<\cdots<t_j<\cdots<t_K<1$ be $K$ points such that  $ {\bf t} =
(t_1,  \ldots, t_K) \in B_{\ep^{3/4}} ( {\bf t}^0)$.
 Set
\begin{equation}
\label{app11}
w_j (x)=  w \left(  \frac{x-t_j}{\ep} \right),
\end{equation}
and
\begin{equation}
\label{r0}
 r_0 =\frac{1}{10} \left( \min \left(t_1^0 +1, 1-t_K^0,
\frac{1}{2}\min_{i \not = j} |t_i^0 -t_j^0|\right)\right).
\end{equation}

Let $\chi: \R \to [0, 1]$ be a smooth cut-off function  such that
$\chi(x)=1$ for $|x|<1$ and $\chi(x)=0$ for $|x|>2$.
We now define the approximate solution as
\begin{equation}
\label{app1}
\tilde{w}_j(x)= w_j (x) \chi\left(\frac{x-t_j}{r_0}\right).
\end{equation}

It is easy to see that $\tilde{w}_j (x)$ satisfies
\begin{equation}
\label{a11}
 \ep^2  \tilde{w}_j^{''} -  \tilde{w}_j + \tilde{w}_j^3 =\mbox{e.s.t.}
\end{equation}
in $L^2(-1,1)$, where e.s.t. denotes an exponentially small term.

Let
\begin{equation}
\label{vector}
\hat{A}=w_{\ep,{\bf t}}(x)=
\sum_{j=1}^K
\frac{1}{\sqrt{v_j^\ep}} \tilde{w_j}(x),\quad \mbox{ where }
v_j^\ep=T[w_{\ep,{\bf t}}](t_j^\ep),
\end{equation}

\begin{equation}
\label{v}
\hat{v}=T[w_{\ep,{\bf t}}],
\end{equation}
where $T[A]$ is defined by
(\ref{ta}) and  ${\bf t} \in B_{\ep^{3/4}} ( {\bf t}^0)$.


Then by (\ref{approx}) we have
\begin{equation}
\label{tauiia}
v_i^\ep:= T[\hat{A}] (t_i^\ep)= \lim_{\ep\to 0} T[w_{\ep,{\bf t}}] (t_i^\ep)
+O\left(\ep\log\frac{1}{\ep}\right).
\end{equation}



\label{tauii}

Now let $x=t_i+\ep y$. We find
for $\hat{A}=w_{\ep, {\bf t}}$:
\[
T[\hat{A}](t_i+\ep y)-T[\hat{A}](t_i)
=
\]\[=\int_{t_i}^{t_i+\ep y}
\frac{1}{ \hat{D} \left(A_0+\frac{1}{\ep}\hat{A}\right)^2} \int_{-L}^s
\left(
\frac{1}{\ep}\hat{v}(\ep A_0+\hat{A})^3-\bar{A}+A_0
\right)
 dt\,ds
\]\[
=\ep^2 \int_{0}^y
\frac{1}{ \hat{D} \left(A_0+\frac{1}{\ep}w(\bar{s})\right)^2} \int_{0}^{\bar{s}}
\frac{1}{\ep}
\frac{(w(\bar{t}))^3}{\sqrt{v_i^\ep}}
 d\bar{t}\,d\bar{s}
 \]\[
+\ep^2 \int_{0}^y
\frac{1}{ \hat{D} \left(A_0+\frac{1}{\ep}w(\bar{s})\right)^2}
\,d\bar{s}
\frac{1}{\ep}
\left[
\int_{-\infty}^{0}
\frac{(w(\bar{t}))^3}{\sqrt{v_i}}\,d\bar{t}
+\sum_{j=1}^{i-1} \frac{\sqrt{2}\pi}{\sqrt{v_j}}H(t_i-t_j)-
\sum_{j=1}^K \frac{\sqrt{2}\pi}{\sqrt{v_j}}
\frac{t_i+L}{2L}
\right]\]\[
\times
\left(1+O\left(\ep\log\frac{1}{\ep}\right)\right)
\]\[
=
\ep P_i^\ep(y)
+\ep
\int_{0}^y
\frac{1}{ \hat{D} \left(A_0+\frac{1}{\ep}w(\bar{s})\right)^2}
\,d\bar{s}
\left[
\frac{1}{2}
\frac{\sqrt{2}\pi}{\sqrt{v_i}}
+\sum_{j=1}^{i-1} \frac{\sqrt{2}\pi}{\sqrt{v_j}}-
\sum_{j=1}^K \frac{\sqrt{2}\pi}{\sqrt{v_j}}
\frac{t_i+L}{2L}
\right]
\]\begin{equation}
\times\left(1+O\left(\ep\log\frac{1}{\ep}\right)\right),
\label{tax}
\end{equation}
in $L^2(\Om_\ep)$,
where
\[
P_i^\ep(y)=\int_{0}^y
\frac{1}{ \hat{D} \left(A_0+\frac{1}{\ep}w(\bar{s})\right)^2}
\int_{0}^{\bar{s}}
\frac{(w(\bar{t}))^3}{\sqrt{v_i^\ep}}
 d\bar{t}\,d\bar{s}
\]
using (\ref{amp1}).
Note that $P_i^\ep$ is an even function and the second term is an odd function in $y$.

We now derive the following estimate for all $y\geq 0$:
\[
(T[w_{\ep, {\bf t}}](t_i+\ep y)-T[w_{\ep, {\bf t}}](t_i))\,w_{\ep, {\bf t}}^3
\leq C
\ep \int_0^y \frac{1}{\hat{D}\left(A_0+\frac{1}{\ep} w(\bar{s})\right)^2}\,d\bar{s}\, w^3(y)
\]\[
\leq \frac{C}{\hat{D}} \ep  \int_0^y \ep^2 \frac{1}{w^2(\bar{s})}\,d\bar{s}\,w^3(y)
\]\[
\leq \frac{C}{\hat{D}} \ep^3 \int_0^y e^{2\bar{s}} \,d\bar{s} \,e^{-3y}
\]\[
\leq\frac{C}{\hat{D}} \ep^3 e^{-y}.
\]
For $y<0$ there is an obvious modification of this estimate. Further, it can be extended to cover both the cases when $w_{\ep, {\bf t}}^3$ is replaced by $\ep w_{\ep, {\bf t}}^2$ or  $\ep^2 w_{\ep, {\bf t}}$, respectively, giving the same upper bound in either case.

Now if we define the following norm
\begin{equation}
\| f\|_{**}= \| f \|_{L^2 (\Omega_\epsilon)} + \sup_{ -  \frac{L}{\epsilon} <y <\frac{L}{\epsilon} } [\max(\min_{i} e^{-\frac{1}{2} |y-\frac{t_i}{\epsilon}| }, \sqrt{\epsilon})]^{-1} | f(y)|
\end{equation}
then by the decay of $w_{\ep, {\bf t}}$ and the definition of the norm, we infer that
\begin{equation}
\|(T[w_{\ep, {\bf t}}](t_i+\ep y)-T[w_{\ep, {\bf t}}](t_i))\,w_{\ep, {\bf t}}^3\|_{**}=
O(\ep^{5/2}),
\label{estw3}
\end{equation}
\begin{equation}
\|(T[w_{\ep, {\bf t}}](t_i+\ep y)-T[w_{\ep, {\bf t}}](t_i))\,\ep w_{\ep, {\bf t}}^2\|_{**}=
O(\ep^{5/2})
\label{estw2},
\end{equation}
\begin{equation}
\|(T[w_{\ep, {\bf t}}](t_i+\ep y)-T[w_{\ep, {\bf t}}](t_i))\,\ep^2 w_{\ep, {\bf t}}\|_{**}=
O(\ep^{5/2}).
\label{estw1}
\end{equation}


Let us now define
\begin{equation}
\label{sa}
S_\ep[\hat{A}] := \ep^{2}\hat{A}_{xx}-\hat{A}+T[\hat{A}] (\varepsilon A_0+\hat{A})^{3}
\end{equation}
where $T[\hat{A}]$ is defined in
(\ref{ta}). Next we set $\hat{A}=w_{\ep, {\bf t}}$ and
compute $S_\ep[w_{\ep, {\bf t}}]$. In fact,
\[S_\ep[w_{\ep, {\bf t}}]=
\ep^{2}
({w_{\ep, {\bf t}}})_{xx}-w_{\ep, {\bf t}}+T[w_{\ep, {\bf t}}]\,(\ep A_0+ w_{\ep, {\bf t}})^{3}
\]\[=\ep^2
({w_{\ep, {\bf t}}})_{xx}-w_{\ep, {\bf t}}+T[w_{\ep, {\bf t}}](t_i)\,w_{\ep, {\bf t}}^{3}
\]\[
+T[w_{\ep, {\bf t}}](\ep^3 A_0^3+3\ep^2 A_0^2 w_{\ep, {\bf t}} + 3 \ep A_0 w_{\ep, {\bf t}}^2)
\]\[
+\left(T[w_{\ep, {\bf t}}](t_i+\ep y)-T[w_{\ep, {\bf t}}](t_i)\right)\, w_{\ep, {\bf t}}^{3}
\]
\begin{equation}
\label{defe}
=:E_1+E_2+E_3.
\end{equation}
We compute
\[
E_1= \sum_{i=1}^K \frac{1}{\sqrt{v_{i}^\ep}}
\left(
\tilde{w}_i^{''}
-\tilde{w}_i+
\frac{T[w_{\ep, {\bf t}}](t_i)}{v_i}
\tilde{w}_{i}^3
\right)
=\mbox{e.s.t.}
\]
in $L_2(\Om_\ep) $ since $v_i^\ep=T[w_{\ep, {\bf t}}](t_i)$.
Further, we get
\[
E_2=O(\ep )
\]
in $L^2(\Om_\ep) $ since $T[w_{\ep, {\bf t}}]$ is bounded in $L^\infty(\Om_\ep)$ and $w_{\ep, {\bf t}}$ is bounded in $L^2(\Om_\ep) $.  We also notice that actually $ E_2= O(\ep e^{-\min_i (|y-\frac{t_i}{\ep}|)})$. Lastly, we derive
\[
E_3=
\sum_{i=1}^K
\frac{1}{(v_{i}^{\ep})^{3/2}}
\left(
T[w_{\ep, {\bf t}}](t_i+\ep y)-
T[w_{\ep, {\bf t}}](t_i)
\right)
\tilde{w}_{\ep, i}^3
=O(\ep^3)
\]
in $L^2(\Om_\ep)$ by (\ref{estw3}).

\medskip

Combining these estimates, we conclude that
\begin{equation}
\label{estsa}
 \|S [ w_{\ep, {\bf t}}]\|_{**} =  O(\ep).
\end{equation}

\begin{remark} {\rm
The estimate (\ref{estsa}) shows that our choice of approximate solution given in (\ref{vector}) and (\ref{v}) is suitable. This will enable us in the next two sections to rigorously construct a steady state which is very close to the approximate solution.
} \end{remark}

\section{Liapunov-Schmidt Reduction}
\setcounter{equation}{0}

 In this section, we use Liapunov-Schmidt reduction
 to solve the  problem
\begin{equation}
\label{ls}
 S_\ep[ w_{\ep, {\bf t}} + v] = \sum_{i=1}^K \beta_i
 \frac{d\tilde{w}_i}{dx}
\end{equation}
for real constants $ \beta_i$ and a function $ v\in H^2(-\frac{1}{\ep},
\frac{1}{\ep})$  which is small in the corresponding norm (to be defined later), where $ \tilde{w}_i$ is given by (\ref{app1}) and $w_{\ep,{\bf t}}$ by
(\ref{vector}). This is the first step in the Liapunov-Schmidt reduction method. We shall follow the general procedure used in \cite{ww}.

To this end, we need to study
the  linearized operator
 \[
\tilde{L}_{\ep, {\bf t}}: H^2 (\Omega_\ep) \to
L^2(\Omega_\ep)
\]
given by
\begin{equation}
\label{lept}
\tilde{L}_{\ep, {\bf t}}
:= S_{\epsilon}^{'} [w_{\ep, {\bf t}}]\phi =
\ep^2\Delta\phi-\phi
+T[w_{\ep, {\bf t}}]
3(\ep A_0+w_{\ep, {\bf t}})^{2}
\phi
+(T'[w_{\ep, {\bf t}}]\phi) (\ep A_0+w_{\ep, {\bf t}})^{3},
\end{equation}
where
for $\hat{A}=w_{\ep, {\bf t}}$ and a given function $\phi\in L^2(\Om)$ we define
 $T^{'}[\hat{A}]\phi$ to be the unique solution of
\begin{equation}
\label{tap}
\left\{
\begin{array}{ll}
\hat{D}
\left(
(A_0 +\frac{1}{\ep} \hat{A})^2
(T'[\hat{A}]\phi)_{x}
\right)_{x}
-\frac{1}{\ep} (T'[\hat{A}]\phi)
(\ep A_0 +\hat{A})^{3}
-\frac{1}{\ep} T[\hat{A}]3
(\ep A_0 +\hat{A})^{2}\phi
=0,\\
\null \hfill\quad\text{for}\quad -L<x<L,
\\[3mm]
(T'[\hat{A}]\phi)_{x}(-L)=(T'[\hat{A}]\phi)_{x}(L)=0.
\end{array}
\right.
\end{equation}

 The norm for the error function $\phi$ is defined as follows

 \begin{equation}
 \label{normdef}
 \| \phi \|_{*}= \| \phi \|_{H^2 (\Omega_\ep)} + \sup_{ -  \frac{L}{\ep} <y <\frac{L}{\ep} }  [\max (\min_{i} e^{-\frac{1}{2}|y-\frac{t_i}{\ep}| }, \sqrt{\ep})]^{-1} | \phi (y)|.
 \end{equation}

We define the approximate kernel and co-kernel, respectively, as follows:
\[ {\cal K}_{\ep, {\bf t}} := \mbox{span} \ \left\{
\frac{d\tilde{w}_i}{dx}
 \Bigg|
i=1,\ldots,K\right\} \subset H^2 (\Om_\ep),
\]
\[ {\cal C}_{\ep, {\bf  t}} := \mbox{span} \ \left\{
\frac{d \tilde{w}_i}{dx}
\Bigg|
i=1,\ldots,K\right\} \subset L^2 (\Om_\ep).
\]
From
(\ref{vecnlep})
we recall the definition of the following system of NLEPs :
\begin{equation}
L\Phi:=\Phi_{yy}-\Phi+3w^2\Phi-3w^3\frac{\int w^2{\cal E}\Phi\,dy}{\int w^3\,dy},
\label{linop1}
\end{equation}
where
\[\Phi=
\left(
\begin{array}{l}
\phi_1\\ \phi_2 \\ \vdots \\ \phi_K
\end{array}
\right)
\in(H^2(\R))^K.
\]
By Lemma 3.3 of \cite{ww}  we know that
\[ L:
 (X_0\oplus  \cdots \oplus X_0)^\perp
\cap (H^2(\R))^K
\to
 (X_0\oplus  \cdots \oplus X_0)^\perp
\cap (L^2(\R))^K\]
is invertible and its inverse is bounded.

We will show that this system is the limit of the operator
$\tilde{L}_{\ep, {\bf t}}$  (defined in (\ref{lept})) as $\ep\to 0$.
We also introduce the projection $\pi_{\ep, {\bf t}}^\perp:
L^2(\Om_\ep)\to {\cal C}^\perp_{\ep, {\bf t}}$
and study the operator $L_{\ep, {\bf t}}:=
\pi_{\ep, {\bf t}}^\perp\circ \tilde{L}_{\ep, {\bf t}}$.
By letting $\ep\to 0$, we will show that $L_{\ep, {\bf t}}:\,
{\cal K}_{\ep, {\bf t}}^\perp \to {\cal C}_{\ep, {\bf t}}^\perp$
is invertible and its inverse is uniformly bounded provided
$\ep$ is small enough.
This statement is contained in the following proposition.
\begin{proposition}
\label{A}
\label{mainprop}
There exist positive constants $\bar{\ep},\,\bar{\delta}, \lambda$
such that for all $\ep\in(0,\bar{\ep})$, ${\bf t} \in\Om^K$ with
$ \min(| 1+t_1|, |1-t_K|, \min_{i \not =j}|t_i-t_j|)>\bar{\delta}$,
\begin{equation}
\label{normesta}
\|L_{\ep, {\bf t} } \phi \|_{**}\geq \lambda
\|\phi\|_{*}.
\end{equation}
Furthermore, the map
\[L_{\ep,{\bf t}}=
\pi_{\ep,{\bf t}}^\perp\circ \tilde{L}_{\ep,{\bf t}}:\,
{\cal K}_{\ep,{\bf t}}^\perp \to {\cal C}_{\ep,{\bf t}}^\perp\]
is surjective.
\end{proposition}

\noindent
{\bf Proof of Proposition \ref{mainprop}:}  This proof uses the method
of Liapunov-Schmidt reduction following for example the approach
in \cite{iww2}, \cite{ww} and \cite{ww13}.

Suppose that (\ref{normesta}) is false. Then there exist
sequences $\{\ep_k\},\,\{{\bf t}^k\},\,\{\phi^k\}$ with $\ep_k\to 0$,
 ${\bf t}^k\in\Om^K$,
 $\min(| 1+t_1^k|, |1-t_K^k|, \min_{i \not =j}|t_i^k-t_j^k|)>\bar{\delta}$,
 $\phi^k=\phi_{\ep_k}\in K_{\ep_k,{\bf t}^k}^\perp$,
 $k=1,2,\ldots$
such that
\begin{eqnarray}
 &&\| L_{\ep_k,{{\bf t}^k }} \phi^k\|_{**}
\to 0 \label{lpt}\qquad\mbox{as }k\to\infty,\\
&&\| \phi^k\|_{*}=1,\abst
k=1,\,2,\,\ldots\,. \label{onnorm}
\end{eqnarray}
We define $\phi_{\ep,i}$, $i=1,2,\ldots,K$ and $\phi_{\ep,K+1}$ as follows:
\begin{equation}
\label{phei}
\phi_{\ep,i}(x)=\phi_\ep(x) \chi\left(\frac{x-t_i}{r_0}\right),\quad x\in\Om,
\end{equation}
\[\phi_{\ep,K+1}(x)=\phi_{\ep}(x)-\sum_{i=1}^K \phi_{\ep,i}(x),\quad x\in\Om.\]
At first (after rescaling) the functions $\phi_{\ep,i}$ are only defined on $\Om_\ep$.
However, by a standard result they can be extended to $\R$ such that their
norm in $H^2(\R)$ is still bounded by a constant independent of $\ep$
and $\bf t$ for
$\ep$ small enough. In the following we will deal with this extension. For
simplicity of notation we keep the same notation for the extension.
Since for $i=1,2,\ldots,K$
each sequence
$\{\phi_i^k\}:=\{\phi_{\ep_k,i}\}$ ($k=1,2,\ldots$)
is bounded in $H^2_{loc}(\R)$ it
converges weakly to a
 limit in $H^2_{loc}(\R)$, and therefore also strongly
in $L^2_{loc}(\R)$ and $L^{\infty}_{loc}(\R)$.
Denoting these limits by $\phi_i$,
then $\phi=\left(\begin{array}{c}\phi_1 \\ \phi_2\\ \vdots \\
\phi_K\end{array}\right)$ solves the system
$L \phi=0.$
By Lemma 3.3 of \cite{ww}, it follows that $\phi\in \mbox{Ker}(L)= X_0\oplus \cdots \oplus X_0$.
Since  $\phi^k\in K_{\ep_k,t_k}^\perp$, taking $k\to\infty$, we get
$\phi\in \mbox{Ker}(L)^\perp$. Therefore, we have $\phi=0$.

By elliptic estimates we get $\|\phi_{\ep_k,i}\|_{H^2(\R)} \to 0$  as
$k\to\infty$
for $i=1,2,\ldots,K$.

Further, $\phi_{\ep,K+1}\to \phi_{K+1}$
in $H^2(\R)$,
where $\Phi_{K+1}$ satisfies
\[ (\phi_{K+1})_{yy}-\phi_{K+1}=0 \quad\mbox{ in }\R.\]
Therefore, we conclude that $\phi_{K+1}=0$ and
$\|\phi_{K+1}^k\|_{H^2(\R)} \to 0$ as $k\to \infty$.

Once we have $ \| \phi_i \|_{H^2 (\R)} \to 0$, the maximum principle implies that $ \|\phi_i \|_{*} \to 0$ since the operator $L_{\ep, {\bf t}}$ essentially behaves like $ \phi_i^{''} - \phi_i$ for $ |x-t_i|>>\ep$. This contradicts the assumption that $\|\phi^k\|_{*}=1$.
To complete the proof of Proposition \ref{A},
 we just need to show that the conjugate operator
to $L_{\ep, {\bf t}}$
(denoted by $ L_{\ep, {\bf t}}^*$) is
injective from ${\cal K}_{\ep, {\bf t}}^\perp $ to ${\cal C}_{\ep,
{\bf t}}^\perp$.
Note that
$
L_{\ep,{\bf t}}^*\phi=\pi_{\ep,{\bf t}}\circ \tilde{L}_{\ep,{\bf t}}^*\phi$
with
\[\tilde{L}_{\ep,{\bf t}}^*\phi:=
\ep^2\Delta\phi-\phi
+T[w_{\ep, {\bf t}}]
3(\ep A_0+w_{\ep, {\bf t}})^{2}
\phi
+(T'[w_{\ep, {\bf t}}]\phi (\ep A_0+w_{\ep, {\bf t}})^{3}).
\]
The  proof for
 $ L_{\ep, {\bf t}}^*$
 follows  the same lines as the proof
for
 $ L_{\ep, {\bf t}}$
and is therefore omitted. Here also the nondegeneracy condition (\ref{condbc}) is required.
For further technical details we refer to \cite{ww}.
\epr

Now we are in the position to solve the equation
\begin{equation}
{\cal \pi}_{\ep, {\bf t}}^\perp \circ S_{\ep} ( w_{\ep, {\bf t}} +
\phi) =0.
\label{solve}
\end{equation}
Since $ L_{\ep, {\bf t}}|_{K_{\ep, {\bf t}}^{\perp}}$ is invertible
(call the inverse $L^{-1}_{\ep, {\bf t}}$) ,we can rewrite this equation as
\begin{equation}
\phi=-(L_{\ep, {\bf t}}^{-1} \circ {\cal \pi}_{\ep,{\bf t}}^\perp
\circ
S_{\ep}
(w_{\ep, {\bf t}}) )
 - (
L_{\ep, {\bf  t}}^{-1}\circ {\cal \pi}_{\ep,
{\bf t}}^\perp\circ N_{\ep, {\bf t}}(\phi))\equiv
M_{\ep, {\bf t}}(\phi),
\label{fix}
\end{equation}
where
\begin{equation}
\label{nep}
N_{\ep,{\bf t}}(\phi)=S_{\ep}(w_{\ep, {\bf t}} +\phi
) -S_{\ep}(w_{\ep, {\bf t}})
-S_{\ep}^{'} (w_{\ep, {\bf t}})
\phi
\end{equation}
and
the operator $M_{\ep, {\bf t}}$ has been defined by (\ref{fix}) for
$\phi\in H^2 (\Omega_{\ep})$. The strategy of the proof is to
show that
the operator $M_{\ep, {\bf t}}$ is a contraction on
\[ B_{\ep,\delta}\equiv\{\phi\in
H^2(\Omega_{\ep})
\,:\, \|\phi\|_{*}<\delta\} \]
if $\ep$ is small enough and $\delta$ is suitably chosen. By
(\ref{estsa}) and Proposition \ref{A} we have that
\[
 \|M_{\ep, {\bf t}}(\phi)\|_{*}
\leq\lambda^{-1}
\left(\|{\cal \pi}_{\ep, {\bf t}}^\perp \circ N_{\ep, {\bf t}}(\phi)
\|_{**} +\left\|{\cal \pi}_{\ep, {\bf t}}^\perp \circ S_\ep
(
w_{\ep, {\bf t}}
)\right\|_{**}\right)
\]\[
\leq \lambda^{-1}C_0(
c(\delta)\delta +\ep),\]
where
$\lambda>0$ is independent of $\delta>0$, $\ep>0$ and
$c(\delta)\to 0$ as $\delta\to 0$. Similarly, we show that
\[
\|M_{\ep,{\bf t}}(\phi)-M_{\ep,{\bf t}}(\phi^{'})\|_{*}
\leq
\lambda^{-1}C_0(
c(\delta)\delta)\|\phi-\phi^{'}\|_{*},\]
where $c(\delta)\to 0$ as $\delta\to 0$. Choosing
$\delta= C_3\ep \mbox{ for $\lambda^{-1}C_0<C_3$}$
and taking $\ep$ small enough,
then  $M_{\ep, {\bf t}}$
maps $B_{\ep,\delta}$ into $B_{\ep,\delta}$ and it is a
contraction mapping in $B_{\ep, \delta}$.
The existence of a fixed point
$\phi_{\ep, {\bf t}}$ now
follows from the standard  contraction mapping principle and
$\phi_{\ep, {\bf t}}$ is a solution of
(\ref{fix}).
\epr

We have thus proved
\begin{lemma}
There exist $\overline{\ep}>0$ $\overline{\delta}>0$
 such that for every pair of $\ep, {\bf t}$
with
$0<\ep<\overline{\ep}$ and ${\bf t}\in\Om^K$,
$1+t_1>\overline{\delta}$, $1-t_K>\overline{\delta}$,
$\frac{1}{2}|t_i-t_j|>
\overline{\delta}$ there is a unique
$\phi_{\ep, {\bf t}}\in K_{\ep, {\bf t}}^{\perp}$
satisfying $S_{\ep}(w_{\ep, {\bf t}} + \phi_{\ep,{\bf t}}
) \in  {\cal C}_{\ep, {\bf t}}$. Furthermore,  the following estimate holds
\begin{equation}
\|\phi_{\ep, {\bf t}}\|_{*}\leq C_3\ep.
\label{estphi}
\end{equation}
\label{lem34}
\end{lemma}


In the next section we  determine the positions of the spikes so that the resulting steady state is an exact solution of the original problem.

\section{The reduced problem}
\setcounter{equation}{0}

In this section we solve the reduced problem and complete the proof of the
existence result for asymmetric spikes in Theorem \ref{existenceas}.

By Lemma \ref{lem34}, for every $ {\bf t} \in B_{\ep^{3/4}} ({\bf t}^0)$,
there exists a unique  $\phi_{\ep, {\bf t} }
\in {\cal K}_{\ep, {\bf t}}^\perp$, solution of
\begin{equation}
\label{see}
S [
w_{\ep, {\bf t}} + \phi_{\ep, {\bf t} }
] =  v_{\ep, {\bf t}}
\in {\cal C}_{\ep, {\bf t}}.
\end{equation}
The idea here is to find  ${\bf t}^\ep=(t_1^\ep,\ldots,t_K^\ep)$ near
${\bf t}^0$ such that also
\begin{equation}
\label{see2}
S  [w_{\ep, {\bf t}^\ep } + \phi_{\ep,  {\bf t}^\ep}
] \perp {\cal C}_{\ep, {\bf t}^\ep}
\end{equation}
(and therefore
$S  [w_{\ep, {\bf t}^\ep } + \phi_{\ep,  {\bf t}^\ep} ]=0$).

 To this end, we let
\[
W_{\epsilon,i}( {\bf t}):=\frac{v_i}{\ep^2}
\int_{-1 }^1
S  [w_{\ep, {\bf t}} +\phi_{\ep, {\bf t}}]
\frac{d \tilde{w}_i}{dx} \, dx
,\]
\[W_{\epsilon}({\bf t}):=(W_{\epsilon,1}({\bf t}),...,
W_{\epsilon,K}({\bf t}) ) : B_{ \ep^{3/4}} ({\bf t}^0) \to \R^K.\]

Then $W_{\epsilon}({\bf t})$ is a  map which is continuous  in ${\bf t}$ and
our problem is reduced to finding
 a zero of the vector field $W_\ep ({\bf t})$.

Let us now calculate $W_{\epsilon}({\bf t})$ as follows:
\[ W_{\ep, i} ( {\bf t})=
\frac{v_i}{\ep^2}
\int_{-1 }^1
S [
w_{\ep, {\bf t}} + \phi_{\ep, {\bf t}}]
\frac{d \tilde{w}_i}{dx}
\]
\[=
\frac{v_i}{\ep^2}
\int_{-1 }^1
S [
w_{\ep, {\bf t}}]
\frac{d \tilde{w}_i}{dx}
\]
\[
+
\frac{v_i}{\ep^2}
\int_{-1 }^1
S_\ep^{'}[w_{\ep, {\bf t}}]
\phi_{\ep, {\bf t}}
\frac{d \tilde{w}_i}{dx}
\]
\[
+
\frac{v_i}{\ep^2}
\int_{-L}^L N_\ep (\phi_{\ep, {\bf t}})
\frac{d \tilde{w}_i}{dx} \]
\[= :I_1 + I_2+ I_3,\]
where $I_1, I_2$ and $I_3$ are defined in an obvious way in the last equality.

\medskip

We will now compute  these three integral terms as $\epsilon \to 0$. The result will be
that $I_1$ is the leading term and $I_2 $ and $ I_3$  are $O\left(\ep\log\frac{1}{\ep}\right)$.

For $I_1$, we have
\[I_1=
\frac{v_i}{\ep^2}
\int_{-L}^L (E_1+E_2+E_3)
\frac{d\tilde{w}_i}{dx}
\, dx  =
\frac{v_i}{\ep^2}
\int_{-L}^L E_3
\frac{d\tilde{w}_i}{dx}
\, dx +O\left(\ep\log\frac{1}{\ep}\right),
\]
where $E_1,\,E_2,\,E_3$ are defined in (\ref{defe}).
For $E_1$ this estimate is obvious. For $E_2$, we use the decomposition
\[
T[w_{\ep, {\bf t}}](\ep^3 A_0^3+3\ep^2 A_0^2 w_{\ep, {\bf t}} + 3 \ep A_0 w_{\ep, {\bf t}}^2)
\]\[
=T[w_{\ep, {\bf t}}]\ep^3 A_0^3
\]\[
+
T[w_{\ep, {\bf t}}](t_i)(3\ep^2 A_0^2 w_{\ep, {\bf t}} + 3 \ep A_0 w_{\ep, {\bf t}}^2)
\]\[
+
(T[w_{\ep, {\bf t}}](t_i+\ep y)- T[w_{\ep, {\bf t}}](t_i))(3\ep^2 A_0^2 w_{\ep, {\bf t}} + 3 \ep A_0 w_{\ep, {\bf t}}^2).
\]
Then we can estimate the first part directly, the second part using the fact that it is an even function in $y$ and the third part using the estimates (\ref{estw2}) and (\ref{estw1}).

From (\ref{estw3}),
we derive
\[
\frac{v_i}{\ep^2}
\int_{-L}^L E_3
\frac{d\tilde{w}_i}{dx}
\,dx\]
\[ =
\frac{v_i}{\ep^2}
\int_{-L/\ep}^{L/\ep} P_i(y)w^3 (y) \frac{w^{'} (y)}{\sqrt{v_i}} \,dy
\]\[
+\frac{v_i}{\ep^2}
\left[
\frac{1}{2}
\frac{\sqrt{2}\pi}{\sqrt{v_i}}
+\sum_{j=1}^{i-1} \frac{\sqrt{2}\pi}{\sqrt{v_j}}-
\sum_{j=1}^K \frac{\sqrt{2}\pi}{\sqrt{v_j}}
\frac{t_i+L}{2L}
\right]
\]\[
\times
\int_{-L/\ep}^{L/\ep} \int_0^y \frac{1}{\hat{D}\left(A_0+\frac{1}{\ep}\hat{A}(\bar{s})\right)^2}\,d\bar{s}
(\ep A_0+\hat{A})^3 \frac{w^{'} (y)}{\sqrt{v_i}}  \chi_i \,dy
\]\[
+O\left(\ep\log\frac{1}{\ep}\right)
\]
\[= \mbox{e.s.t.}-
 \frac{1}{\hat{D}}
 \left[
\frac{1}{2}
\frac{\sqrt{2}\pi}{\sqrt{v_i}}
+\sum_{j=1}^{i-1} \frac{\sqrt{2}\pi}{\sqrt{v_j}}-
\sum_{j=1}^K \frac{\sqrt{2}\pi}{\sqrt{v_j}}
\frac{t_i+L}{2L}
\right]
+O\left(\ep\log\frac{1}{\ep}\right),
\]
where $\chi_i(x)=\chi\left(\frac{x-t_i}{r_0}\right)$.
Here we have used the fact that $P_i(y)$ is an even function and have computed the following integral
\[
\int_{-L/\ep}^{L/\ep} \int_0^y \frac{1}{\hat{D}\left(A_0+\frac{1}{\ep}\hat{A}(t_i+\ep\bar{s})\right)^2}\,d\bar{s}
(\ep A_0+\hat{A})^3 \frac{w^{'} (y)}{\sqrt{v_i}}  \chi_i(t_i+\ep y) \,dy
\]\[
=\int_{-L/\ep}^{L/\ep} \int_0^y \frac{1}{\hat{D}\left(A_0+\frac{1}{\ep}\hat{A}(t_i+\ep \bar{s})\right)^2}\,d\bar{s}
\frac{1}{4}\frac{d}{dy}(\ep A_0+\hat{A})^4 \chi_i(t_i+\ep y)  \,dy+O\left(\ep^3\log\frac{1}{\ep}\right)
\]\[
=-\frac{\ep^2}{\hat{D}} \int_{-L/\ep}^{L/\ep}
\frac{1}{4}\left(\ep A_0+\hat{A}(t_i+\ep y)\right)^2
\,dy+O\left(\ep^3\log\frac{1}{\ep}\right)
\]\[
=-\frac{\ep ^2}{\hat{D}}\int_{\R}\frac{(w(y))^2}{4v_i}\,dy+O\left(\ep^3\log\frac{1}{\ep}\right)
\]\[
=-\frac{\ep^2}{\hat{D}v_i} \int_{\R} \frac{1}{2\cosh^2 y}\,dy +O\left(\ep^3\log\frac{1}{\ep}\right)
\]\[
=-\frac{\ep^2}{ \hat{D}v_i}+O\left(\ep^3\log\frac{1}{\ep}\right).
 \]


In summary, we have
\begin{equation}
\label{esti1}
I_1=-\frac{1}{\hat{D}}
 \left[
\frac{1}{2}
\frac{\sqrt{2}\pi}{\sqrt{v_i}}
+\sum_{j=1}^{i-1} \frac{\sqrt{2}\pi}{\sqrt{v_j}}-
\sum_{j=1}^K \frac{\sqrt{2}\pi}{\sqrt{v_j}}
\frac{t_i+L}{2L}
\right]
+O\left(\ep\log\frac{1}{\ep}\right).
\end{equation}

For $I_2$, we calculate
\[ I_2=
\frac{v_i}{\ep^2}
\int_{-1 }^1
S^{'}[w_{\ep, {\bf t}}] (\phi_{\ep, {\bf t}})
\frac{d\tilde{w}_i}{dx}
\]
\[=
\frac{v_i}{\ep^2}
\int_{-L}^L
\left[
\ep^2\Delta\phi_{\ep, {\bf t}}-\phi_{\ep, {\bf t}}
+T[w_{\ep, {\bf t}}]
3(\ep A_0+w_{\ep, {\bf t}})^{2}
\phi_{\ep, {\bf t}}
+(T'[w_{\ep, {\bf t}}]\phi_{\ep, {\bf t}}) (\ep A_0+w_{\ep, {\bf t}})^{3}
\right]
\frac{d\tilde{w}_i}{dx}
\]\[
=
\frac{v_i}{\ep^2}
\int_{-L}^L
\left[
\ep^2\Delta
\frac{d\tilde{w}_i}{dx}
-
\frac{d\tilde{w}_i}{dx}
+
3\tilde{w}_{i}^{2}
\frac{T[w_{\ep, {\bf t}}](t_i)}{v_i}d\frac{\tilde{w}_i}{dx}
\right]
\phi_{\ep, {\bf t}}\,dx
\]\[
+
\frac{v_i}{\ep^2}
\int_{-L}^L
\frac{T[w_{\ep, {\bf t}}](x)-T[w_{\ep, {\bf t}}](t_i)}{v_{i}}
\,3\tilde{w}_{i}^{2}
\phi_{\ep, {\bf t}}
\frac{d\tilde{w}_i}{dx}
\,dx
\]\[
+\frac{v_i}{\ep^2}
\int_{-L}^L
T[w_{\ep, {\bf t}}](t_i)3\left(\ep^2 A_0^2+2\ep A_0\frac{\tilde{w}_i}{\sqrt{v_i}}\right)
\phi_{\ep, {\bf t}}
\frac{d\tilde{w}_i}{dx}
\,dx
\]\[
+\frac{v_i}{\ep^2}
\int_{-L}^L
(T[w_{\ep, {\bf t}}](x)-T[w_{\ep, {\bf t}}](t_i))3\left(\ep^2 A_0^2+2\ep A_0\frac{\tilde{w}_i}{\sqrt{v_i}}\right)
\phi_{\ep, {\bf t}}
\frac{d\tilde{w}_i}{dx}
\]\[
+
\frac{v_i}{\ep^2}
\int_{-L}^L
(T'[w_{\ep, {\bf t}}]\phi_{\ep, {\bf t}})
\ep^3A_0^3
\frac{d\tilde{w}_i}{dx}
\,dx
\]\[
+
\frac{v_i}{\ep^2}
\int_{-L}^L
(T'[w_{\ep, {\bf t}}]\phi_{\ep, {\bf t}})(t_i)
 \left(
 3\ep^2 A_0^2 \frac{\tilde{w}_i}{\sqrt{v_i}}+3 \ep A_0 \frac{\tilde{w}_{i}^{2}}{v_i}+\frac{\tilde{w}_{i}^{2}}{v_i^{3/2}}
 \right)
\frac{d\tilde{w}_i}{dx}
\,dx
\]
\[
+
\frac{v_i}{\ep^2}
\int_{-L}^L
[
(T'[w_{\ep, {\bf t}}]\phi_{\ep, {\bf t}})(x)
-
(T'[w_{\ep, {\bf t}}]\phi_{\ep, {\bf t}})(t_i)
]
\]\[
\times \left(
3\ep^2 A_0^3 \frac{\tilde{w}_i}{\sqrt{v_i}}+3 \ep A_0 \frac{\tilde{w}_{i}^{2}}{v_i}+\frac{\tilde{w}_{i}^{2}}{v_i^{3/2}}
\right)
\frac{d\tilde{w}_i}{dx}
\,dx
\]\[
=I_2^1+I_2^2+I_2^3+I_2^4+I_2^5+I_2^6+I_2^7.
\]

\medskip

With obvious notations, we now show that each one of the seven terms is  $O\left(\ep\log\frac{1}{\ep}\right)$ as $\epsilon \to 0$.

For $I_2^1$, it follows from $T[w_{\ep, {\bf t}}]=v_i$, while for $I_2^2$, we use (\ref{estw3}) and the fact that
\begin{equation}
\left\|\frac{\phi_{\ep,{\bf t}}}{w_{\ep,{\bf t}}}\right\|_{L^\infty(\Om_\ep)}\leq C\ep.
\label{estphiw}
\end{equation}
For $I_2^3$, we use $\|T[w_{\ep,{\bf t}}]\|_{L^\infty(\Om_\ep)}=O(1)$ and the fact that $\phi_{\ep,{\bf t}}$ is an even function.
For $I_2^4$, we use (\ref{estw2}), (\ref{estw1}) and (\ref{estphiw}).
For $I_2^5$, the estimate is derived from  $\|T'[w_{\ep,{\bf t}}]\phi_{\ep, {\bf t}}\|_{L^\infty(\Om_\ep)}=O(\ep)$.
For $I_2^6$, we use $(T'[w_{\ep,{\bf t}}]\phi_{\ep, {\bf t}})(t_i)=O(\ep)$ and the fact that
$\tilde{w}_i$ is even.
Lastly, for  $I_2^7$, we use estimates similar to (\ref{estw3}), (\ref{estw2}), (\ref{estw1}) with $T'[w_{\ep,{\bf t}}]\phi_{\ep, {\bf t}}$ instead of  $T[w_{\ep,{\bf t}}]$ and the inequality
\[
\left\|\frac{T'[w_{\ep,{\bf t}}]\phi_{\ep, {\bf t}}}{T[w_{\ep,{\bf t}}]}\right\|_{L^\infty(\Om_\ep)}\leq C\ep.
\]

By arguments similar to the ones for $I_2$, we derive
\begin{equation}
\label{i33}
I_3 = O\left(\ep\log\frac{1}{\ep}\right) \quad\mbox{ in } L^2(\Om_\ep).
\end{equation}
Combining the estimates for $I_1$, $I_2$ and $I_3$, we have
\[ W_{\ep,i} ({\bf t})=
-\frac{1}{\hat{D}}
 \left[
\frac{1}{2}
\frac{\sqrt{2}\pi}{\sqrt{v_i}}
+\sum_{j=1}^{i-1} \frac{\sqrt{2}\pi}{\sqrt{v_j}}-
\sum_{j=1}^K \frac{\sqrt{2}\pi}{\sqrt{v_j}}
\frac{t_i+L}{2L}
\right]
+O\left(\ep\log\frac{1}{\ep}\right).
\]
\begin{equation}
\label{formulaf}
= -\frac{1}{\hat{D}} F_i ({\bf t}) + O\left(\ep\log\frac{1}{\ep}\right),
\end{equation}
where  $F_i ({\bf t})$  was defined in (\ref{fi}).

By assumption (\ref{amp3}),
we have  $  F ({\bf t}^0) =0$.
Next we show that
\[ \det ( \nabla_{{\bf t}^0} F ({\bf t}^0)) \not = 0\]
in the case of two spikes with amplitudes $v_1=v_s<v_l=v_2$. We compute
\[ \nabla_{{\bf t}} F ({\bf t})=
D_{\bf t} F+(D_{v} F) (D_{\bf t} v),
\]
where
\[
D_{\bf t} F
=-\left(\sum_{j=1}^2 \frac{\sqrt{2}\pi}{\sqrt{v_j}}\frac{1}{2L}\right){\cal I},
\]\[
D_{v} F=\frac{\sqrt{2}\pi}{4}
\left(
\begin{array}{cc}
0 & v_2^{-3/2} \\[3mm]
-v_1^{-3/2} & 0
\end{array}
\right),
\]\[
D_{\bf t} v= \frac{\pi}{\sqrt{2}\hat{D} A_0^2}
\,
\frac{1}{\left(\left(\frac{v_2}{v_1}\right)^{3/2}+1-\frac{\pi L}{\sqrt{2}\hat{D}A_0^2 v_1^{3/2}}\right)}
\left(
\begin{array}{cc}
v_1^{-1/2} & v_2^{-1/2} \\[3mm]
-v_1^{-2} v_2^{3/2} & -v_1^{-3/2} v_2
\end{array}
\right).
\]
This implies, using (\ref{amp1}) and (\ref{amp6}), that
\[\nabla_{{\bf t}^0} F ({\bf t}^0)=
\frac{\pi^2}{4 \hat{D} A_0^2}\,
 \frac{1}{\left(\left(\frac{v_2}{v_1}\right)^{3/2}+1-\frac{\pi L}{\sqrt{2}\hat{D}A_0^2 v_1^{3/2}} \right)}
\]\[
\times
\left(
\begin{array}{cc}
-v_1^{-5/2}v_2^{1/2} -v_1^{-1}v_2^{-1}+v_1^{-2}+2v_1^{-3/2}v_2^{-1/2} & -v_1^{-3/2}v_2^{-1/2}   \\[3mm]
-v_1^{-2} & -v_1^{-5/2}v_2^{1/2} -v_1^{-1}v_2^{-1}+2v_1^{-2}+v_1^{-3/2}v_2^{-1/2}
\end{array}
\right).
\]
Next, we compute
\[\det ( \nabla_{{\bf t}^0} F ({\bf t}^0))=
\frac{\pi^4}{16 (\hat{D} A_0^2)^2}\,\frac{1}{v_1^2v_2^2}\,
 \frac{1}{\left(\alpha^3-2\alpha^2-2\alpha+1 \right)^2}
 \]\[
 \times
\left(
(\alpha^3-\alpha^2-2\alpha+1)(\alpha^3-2\alpha^2-\alpha+1)-\alpha^3
\right)
\]
\[=
\frac{\pi^4}{16 (\hat{D} A_0^2)^2}\,
\frac{1}{v_1^2v_2^2}
\,\frac{\alpha^3-\alpha^2-\alpha+1}{\alpha^3-2\alpha^2-2\alpha+1},
\]
where $\alpha=\sqrt{\frac{v_2}{v_1}}$.
Therefore, we have $\det ( \nabla_{{\bf t}^0} F ({\bf t}^0))\not = 0$, except  for two specific positive values of $\alpha$:
$\alpha=1$ (the bifurcation point of asymmetric from symmetric spikes which is not included in Theorem \ref{existenceas})
and $\alpha=\frac{1+\sqrt{5}}{2}$ (corresponding to the eigenvalue $e_{m,1}=\frac{3}{2}$ in Section 5 which has been excluded from Theorem \ref{existenceas}).

Thus, under the conditions of Theorem \ref{existenceas},  we get \[W_\ep ({\bf t})
= -\frac{1}{\hat{D}} \nabla_{{\bf t}^0} F ({\bf t}^0)({\bf t}-{\bf t}^0) + O\left(|{\bf t} -{\bf t}^0|^2+\ep\log\frac{1}{\ep}\right).\]

Since $W_\ep({\bf t})$ is continuous in $\bf t $,
standard degree theory \cite{danclec} implies
that for
$\epsilon$ small enough and $\delta$ suitable chosen there exist
${\bf t^{\epsilon}}\in B_\delta({\bf t}^0)$
such that
$W_{\epsilon}({\bf t^{\epsilon}})=0$ and ${\bf  t^\ep} \to {{\bf t}^0}$.
For further technical details of the argument, we refer to \cite{ww13}.

\epr

Thus we have proved the following proposition.
\begin{proposition}
\label{redprob}
For $\epsilon$ small enough, there exist  points
${\bf t}^{\epsilon}$ with
${ \bf t}^{\epsilon}\to { \bf t}^0$ such that
$W_{\epsilon}({\bf t}^\ep)=0$.
\label{conver}
\end{proposition}


Finally, we complete the proof of Theorems \ref{existencesym} and \ref{existenceas}.

\noindent
{\bf Proof of Theorem \ref{existenceas}:} By Proposition \ref{redprob},
there exist ${\bf t}^\ep \to {\bf t}^0$
such that $ W_\ep ({\bf t}^\ep) =0$.
In other words, $S [w_{\ep, {\bf t}^\ep}+\phi_{\ep,
{\bf t}^\ep}] =0$. Let $\hat{A}_\ep =
w_{\ep, {\bf t}^\ep} +\phi_{\ep, {\bf t}^\ep},\, \hat{v}_\ep = T[w_{\ep, {\bf t}^\ep} + \phi_{\ep, {\bf t}^\ep} ]$.
By  the maximum principle, we conclude that $\hat{A}_{\ep} >0,\, \hat{v}_\ep >0$.
Moreover
$(\hat{A}_\ep, \hat{v}_\ep)$ satisfies all the properties of Theorem \ref{existenceas}.

\epr

\noindent
{\bf Proof of Theorem \ref{existencesym}:}
To prove Theorem \ref{existencesym}, we first construct a single spike in the interval $\left(-\frac{L}{K},\frac{L}{K}\right)$ as above. Then we continue the single spike periodically to a function in the interval $(-L,L)$ and get a symmetric multiple spike in the interval $(-L,L)$.

\epr

\section{Proof of Theorem \ref{existenceani}}
\setcounter{equation}{0}

The proof of Theorem \ref{existenceani} goes exactly as that of Theorem \ref{existenceas}.

\medskip

First, let us derive the location of the single spikes formally: in the first equation in (\ref{sys1}) the term $\varepsilon^3 A_0^{''}$ is very small and can be omitted in the computations. Thus we may assume that
\begin{equation}
\label{n100}
\hat{A} \sim \xi^{-1/2} w\left(\frac{x-t_0}{\varepsilon}\right), \  v(t^0)= \xi
\end{equation}

\medskip

Substituting the above expressions into the second equation of (\ref{sys1}) and noting that $ \hat{v} \frac{1}{\varepsilon} (\varepsilon A_0 +\hat{A})^3 \sim \xi^{-1/2} (\int w^3) \delta_{t_0}$, we see that $\hat{v}$ satisfies in leading order
\begin{equation}
\hat{D} ( A_0^2 \hat{v}_x)_{x} - \xi^{-1/2} \left(\int_\R w^3\,dy\right) \delta_{t_0} + \gamma (x)=0.
\end{equation}
Solving the above equation, we then obtain
\begin{equation}
\label{eqn100}
\hat{v}_x (t_0 -)= - \frac{1}{\hat{D} (A_0 (t_0))^2} \int_{-L}^{t_0} \gamma (x) dx, \ \ \hat{v}_x (t_0 +)=  \frac{1}{\hat{D} (A_0 (t_0))^2} \int_{t_0}^{L} \gamma (x) dx.
\end{equation}

Substituting (\ref{n100}) into the first equation of (\ref{sys1}) and rescaling $x= t_0+\varepsilon y$, we deduce that the error becomes
\begin{equation}
\varepsilon ( \hat{v}_x (t_0-) y^{-} + \hat{v}_x (t_0+) y^{+})  w^3 (y) + \varepsilon A_0(t_0) 3 w^2
 +O(\varepsilon^2),
 \end{equation}
where $y^{-}=\min(y,0)$ and $y^{+}=\max(y,0)$.
 from which we conclude that a necessary condition for the existence of a spike at $t_0$ is that
 \begin{equation}
\int_\R \left[\varepsilon ( \hat{v}_x (t_0-) y^{-} + \hat{v}_x (t_0+) y^{+})  w^3 (y) + \varepsilon A_0(t_0) 3w^2 +O(\varepsilon^2)\right] w^{'} (y)\,dy=0
 \end{equation}
 whence
 \begin{equation}
 \hat{v}_x (t_0-)+\hat{v}_x (t_0+)=0
 \end{equation}
 which is equivalent to (\ref{tep000}). It turns out that (\ref{tep000}) is also sufficient, since the derivative of $\int_{-L}^{t_0} \gamma (x) dx -\int_{t_0}^L \gamma (x) dx$ with respect to $t_0$ is $\gamma (t_0)$ which is strictly positive. The rest of the proof goes exactly as in the proof of Theorem \ref{existenceas}. We omit the details.

\section{Discussion}

In this article we have provided a rigorous mathematical analysis of the formation of spikes in the model of Short, Bertozzi and Brantingham \cite{sbb}. Thus, we have shown that this model naturally leads to the formation of criminality hot-spots. The existence of such hotspots is one of the main stylized facts about criminality. It is observed for an array of criminal activity types. Hot-spots are extensively reported and discussed  in the criminology literature. We refer for example to the articles \cite{je} and \cite{bb} as well as to the references therein.
Now, the fact that a mathematical model yields such hotspots can be viewed as passing one benchmark of validity. The findings in our paper provide such a test for the Short, Bertozzi and Brantingham model \cite{sbb}.

Furthermore, the rigorous analysis carried here sheds light on the mechanisms for the formation of hotspots in this model and the way it quantitatively depends on the parameters.
This type of analysis can then be applied to study issues such as the reduction of hotspots by crime prevention strategies or optimal use of resources to this effect.
One of the goals is to understand when policing strategies actually reduce criminal activity and when they merely displace hot-spots to new areas.

In this paper we have proved  three main new results. First, we showed that we can reduce the quasilinear chemotaxis problems to a Schnakenberg type reaction-diffusion system and derived the existence of symmetric $k$ spikes. Next,  we established the existence of asymmetric spikes in the isotropic case. Lastly, we have studied the pinning effect by inhomogeneous media $A_0 (x)$ and $\gamma (x)$. The stability of these spikes is an interesting open issue.

In \cite{kww1} spikes in two space dimensions are considered by formal matched asymptotics. Our approach of rigorous justification can be extended to that case in a radially symmetric setting, i.e. if the domain is a disk and we construct a single spike located at the centre. We remark that in \cite{kww1} it is assumed that in the outer region (away from the spikes) the system in leading order is semi-linear which allows an extension of the results for the Schnakenberg model to this case. In \cite{kww1}, for the inner region, a numerical computation by solving a core problem yields the profile of the spike.

An alternative approach to the problem in one space dimension would be to write it as a first-order semi-linear ODE system and then apply standard methods, e.g. dynamical system methods for the problem on the real line. This approach becomes cumbersome when we impose Neumann boundary conditions and we add inhomogeneity. We nevertheless refer to a recent paper \cite{HDK}, where the dynamical systems approach is used to construct traveling wave solutions of a quasi-linear  reaction-diffusion-mechanical system.

We remark that there are very few results concerning the analysis of spikes in   quasi-linear reaction diffusion systems. As far as we know, there are two such types of systems. The first one is the chemotaxis system of Keller-Segel type. We refer to \cite{HP} for the background of chemotaxis models and \cite{KW2} for the analysis of spikes to these systems. The other one is the Shigesada-Kawasaki-Teramoto model of species segregation (\cite{SKT}). For the analysis of spikes in a cross-diffusion system, we refer to  \cite{KW1}, \cite{LN} and \cite{WX}.

A family of related models for the diffusion of criminality has been proposed in \cite{bn}. We analyze the formation of hot spots in this class of models in our forthcoming work \cite{bw}. The equations in \cite{bn} also envision the possibility of non-local diffusion. Indeed, social influence can be exercised at long range and it is natural to consider descriptions that take long range diffusion into account. Such a non-local system arising in \cite{bn} reads:
\begin{equation}\label{crime_non-local}
\begin{cases}
 &s_t(x,t)  ={\mathcal L} s(x,t) -s(x,t)+s_b+\alpha(x) u(x,t)\\
 &u_t(x,t) = \Lambda(s) - u(x,t).
\end{cases}
\end{equation}
The case when ${\mathcal L}=\Delta $ is a local diffusion operator provides the framework of the study in
\cite{brr}.  Here, ${\mathcal L}$ can also be a non-local operator such as the fractional Laplace operator or a general non-local interaction term:
$$
{\mathcal L} s (x,t)= \int J(x,y) (s(y,t) - s(x,t) ) dy.
$$
Observe that the steady states reduce to a single non-local equation:

\begin{equation}\label{eqs-stat}
- {\mathcal L} s =  s_b (x)  -s +\alpha(x) \Lambda (s).
\end{equation}

We note that the interaction between non-local diffusion and the mechanism for the formation of spikes is completely open. In particular, the description of the formation of spikes in (\ref{crime_non-local}) and (\ref{eqs-stat}) are open problems. We expect that the decay of the kernel may come into play for the formation of spikes.

\bigskip

\noindent {\bf Acknowledgment.}
The research of Henri Berestycki leading to these results has received funding from the European Research Council under the European Union's Seventh Framework Programme (FP/2007-2013) / ERC Grant Agreement
n.321186 - ReaDi -Reaction-Diffusion
Equations, Propagation and Modelling. Part of this work was done while Henri Berestycki was visiting the Department of Mathematics at University of Chicago. He was also supported by an NSF FRG grant DMS - 1065979, "Emerging issues in the Sciences Involving Non-standard Diffusion".
Juncheng Wei was supported by a GRF grant from RGC of Hong Kong and a NSERC Grant from Canada. Matthias Winter thanks the Department of Mathematics of The Chinese University of Hong Kong for its kind hospitality. Lastly, the authors are thankful to the referees for careful reading of the manuscript and many constructive  suggestions.

\end{document}